\documentclass[11pt]{article}
\usepackage{amsmath}
\usepackage{dsfont}
\usepackage{amsmath}
\usepackage{amsmath,bm}
\usepackage{mathrsfs}
\usepackage{amsmath,amssymb}
\usepackage{amsfonts}
\usepackage{hyperref}
\usepackage{amsthm}
\usepackage{graphicx}
\usepackage{subfigure}
\usepackage{xcolor}
\usepackage{appendix}
\usepackage{cite}

\usepackage{tikz}

\usetikzlibrary{calc,shapes,arrows}

\usepackage{pgfplots}

\usetikzlibrary{patterns}

\usepgfplotslibrary{dateplot}

\renewcommand{\qed}{\hfill\small{$\square$}\normalsize}

\hfuzz=\maxdimen
\theoremstyle{definition}
\newtheorem{lemma}{Lemma}[section]
\newtheorem{definition}{Definition}
\newtheorem{proposition}[lemma]{Proposition}
\newtheorem{theorem}[lemma]{Theorem}
\newtheorem{corollary}[lemma]{Corollary}

\newtheorem{remark}{Remark}
\newtheorem{question}{Question}

\numberwithin{equation}{section}
\renewcommand{\proof}{\textbf{Proof. }}
\renewcommand{\qed}{\hfill\small{$\square$}\normalsize}

\DeclareFixedFont{\Acknowledgment}{OT1}{cmr}{bx}{n}{14pt}
\begin{document}

\title{\bf A new class of discrete conformal structures on surfaces with boundary}
\author{Xu Xu}
\maketitle

\begin{abstract}
We introduce a new class of discrete conformal structures on surfaces with boundary, which
have nice interpolations in 3-dimensional hyperbolic geometry.
Then we prove the global rigidity of the new discrete conformal structures using variational principles,
which is a complement of Guo-Luo's rigidity of the discrete conformal structures in \cite{GL}
and Guo's rigidity of vertex scaling in \cite{Guo2} on surface with boundary.
As a result, new convexities of the volume of generalized hyperbolic pyramids with
right-angled hyperbolic hexagonal bases are obtained.
Motivated by Chow-Luo's combinatorial Ricci flow and Luo's combinatorial Yamabe flow on closed surfaces,
we further introduce combinatorial Ricci flow and combinatorial Calabi flows
to deform the new discrete conformal structures on surfaces with boundary.
The basic properties of these combinatorial curvature flows are established.
These combinatorial curvature flows provide effective algorithms
for constructing hyperbolic metrics on surfaces with totally geodesic boundary components of prescribed lengths.
\end{abstract}




\section{Introduction}

Discrete conformal structure on polyhedral manifolds is a discrete analogue of the well-known conformal structure
on smooth Riemannian manifolds,
which assigns the discrete metrics defined on the edges by scalar functions defined on the vertices.
Since the famous work of William Thurston on circle packings on closed surfaces \cite{T},
different types of discrete conformal structures on closed surfaces has been extensively studied.
See, for instance, \cite{BPS, CL, T, GGLSW, GLSW, GLW, Guo1,  L1, L3, LW, GT,  SWGL, Wu, WGS, WZ, X1, X3, S} and others.

However, the discrete conformal structures on surfaces with boundary are seldom studied.
Motivated by Thurston's circle packings on closed surfaces,
Guo-Luo \cite{GL} first introduced some generalized circle packing type
hyperbolic discrete conformal structures on surfaces with boundary.
Following Luo's vertex scaling of piecewise linear metrics on closed surfaces \cite{L1},
Guo \cite{Guo2} introduced a class of hyperbolic discrete conformal structures, also called vertex scaling,
 on surfaces with boundary.
In this paper, we introduce a new class of hyperbolic discrete conformal structures on ideally triangulated surfaces with boundary,
which has nice geometric interpolations in 3-dimensional hyperbolic geometry.
Then we study the rigidity and deformation of the new discrete conformal structures on surfaces with boundary.

\begin{figure}[!htb]
\centering
  \includegraphics[height=.52\textwidth,width=.5\textwidth]{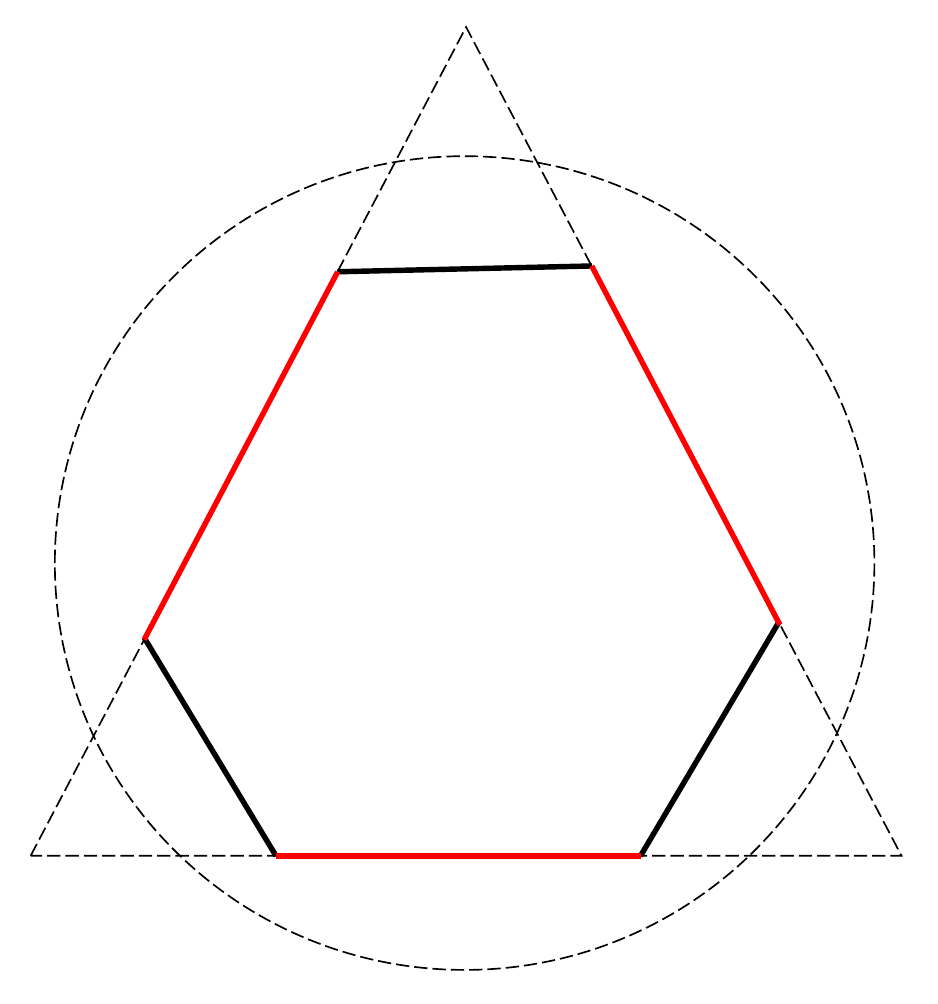}
  \caption{Colored right-angled hyperbolic  hexagon}
  \label{Colored right-angled hyperbolic  hexagon}
\end{figure}

Suppose $\Sigma$ is a compact surface with boundary $B$, which is composed of $n$ boundary components.
$\mathcal{T}$ is an ideal triangulation of $\Sigma$, which can be constructed as follows.
Suppose we have a finite disjoint union of colored topological hexagons, three non-adjacent edges of each hexagon are colored red
and the other edges are colored black.
Please refer to Figure \ref{Colored right-angled hyperbolic  hexagon} for a colored hexagon.
 Identifying the red edges of colored hexagons in pairs by homeomorphisms
gives rise to a quotient space, called an ideal triangulated compact surface with boundary.
The image of each colored hexagon is a face in the triangulation $\mathcal{T}$ and
the image of each red edge in the colored hexagon is an edge in the triangulation $\mathcal{T}$.
The image of the black edges are referred as boundary arcs.
For simplicity, we denote  the boundary components of $(\Sigma, \mathcal{T})$ as $B=\{1,2, \cdots, n\}$,
denote the set of edges in $(\Sigma, \mathcal{T})$ as $E$ and denote the set of faces in $(\Sigma, \mathcal{T})$ as $F$.
An edge connecting the boundary components  $i, j\in B$ is denoted by $\{ij\}$ and
a face adjacent to $i,j,k\in B$ is denoted by $\{ijk\}$.

A basic fact from hyperbolic geometry \cite{Rat} is that given any three positive numbers,
there exists a unique right-angled hyperbolic  hexagon up to hyperbolic isometry
with the lengths of three non-adjacent edges in the hexagon given by the three positive numbers.
Therefore, if $l: E\rightarrow (0,+\infty)$ is a positive function defined on $E$,
every face in $F$ can be realized as a unique right-angled hyperbolic  hexagon up to isometry with the lengths of edges in $E$ given by $l$.
By gluing the right-angled hyperbolic  hexagons along the edges in $E$ in pairs by isomorphisms according to
the ideal triangulation $\mathcal{T}$, we get a hyperbolic metric on the ideally triangulated surface $(\Sigma, \mathcal{T})$
with totally geodesic boundary components.
Conversely, every hyperbolic metric on an ideally triangulated surface $(\Sigma, \mathcal{T})$ with totally geodesic boundary components
with $\mathcal{T}$ geometric
determines a unique map $l: E\rightarrow (0,+\infty)$ with $l_{ij}$ given by the length of the shortest geodesic connecting
the boundary components $i, j\in B$.
The map $l: E\rightarrow (0,+\infty)$ is called as a\textbf{ discrete hyperbolic metric} on $(\Sigma, \mathcal{T})$.
The length $K_i$ of the boundary component $i\in B$ is called
the \textbf{generalized combinatorial curvature} of the discrete hyperbolic metric $l: E\rightarrow (0,+\infty)$ at $i\in B$.

Note that every colored  right-angled hyperbolic  hexagon in the hyperbolic space corresponds to a unique
generalized hyperbolic triangle in the extended hyperbolic space (using the Klein model)
with three hyper-ideal vertices and the three segments between
the hyper-ideal vertices intersecting with the hyperbolic plane $\mathbb{H}^2$.
Please refer to Figure \ref{Colored right-angled hyperbolic  hexagon}.
Further note that for the colored right-angled hyperbolic  hexagons adjacent to the same boundary components,
the corresponding generalized triangles are adjacent to same hyper-ideal vertex.
In this sense, the ideally triangulated hyperbolic surface with boundary $(\Sigma, \mathcal{T}, l)$
can be taken as a triangulated closed surface in the extended hyperbolic space,
with the totally geodesic boundary components corresponding to the hyper-ideal vertices.
For simplicity, a generalized hyperbolic triangle is always referred to a right-angled hyperbolic hexagon in the following,
if it causes no confusion in the context.
Recall that a discrete conformal structure on a triangulated closed surface assigns
the discrete metrics defined on the edges by functions defined on the vertices.
Motivated by the following hyperbolic discrete conformal structures
introduced by Glickenstein-Thomas \cite{GT}, Zhang-Guo-Zeng-Luo-Yau-Gu \cite{ZGZLYG} and Bobenko-Pinkall-Springborn \cite{BPS} on triangulated closed surfaces
 \begin{equation*}
  \begin{aligned}
  l_{ij}=\cosh^{-1}\left(\sqrt{(1+\varepsilon_ie^{2u_i})(1+\varepsilon_je^{2u_j})}+\eta_{ij}e^{u_i+u_j}\right)
  \end{aligned}
  \end{equation*}
with $\varepsilon_i,\varepsilon_j\in \{-1,0,1\}$ and $\eta_{ij}\in \mathbb{R}$,
we introduce the following discrete conformal structures on ideally triangulated surfaces with boundary.

\begin{definition}\label{defn of DCG on surf w bdy}
Suppose $(\Sigma, \mathcal{T})$ is an ideally triangulated surface with boundary
and $\eta: E\rightarrow (-1, +\infty)$ is a weight defined on the edges.
A discrete conformal structure on $(\Sigma, \mathcal{T})$ is a function $u:B\rightarrow \mathbb{R}$
  such that
\begin{equation}\label{DCS on bordered surface}
\begin{aligned}
 l_{ij}
=\cosh^{-1}\left(e^{u_i+u_j}+\eta_{ij}\sqrt{(1+e^{2u_i})(1+e^{2u_j})}\right)
\end{aligned}
\end{equation}
determines a discrete hyperbolic metric $l: E\rightarrow (0, +\infty)$ on $(\Sigma, \mathcal{T})$.
The function $u:B\rightarrow \mathbb{R}$ is called a discrete conformal factor.
\end{definition}


A basic problem in discrete conformal geometry is to understand
the relationships between the discrete conformal structures and their combinatorial curvatures.
We prove the following result on the rigidity of
the discrete conformal structures in Definition \ref{defn of DCG on surf w bdy}.

\begin{theorem}\label{main thm rigidity}
Suppose $(\Sigma, \mathcal{T})$ is an ideally triangulated surface with boundary
and $\eta: E\rightarrow (-1, +\infty)$ is a weight defined on the edges.
If the weight $\eta$ satisfies the following structure condition
  \begin{equation}\label{structure condition}
\begin{aligned}
\eta_{ij}+\eta_{ik}\eta_{jk}\geq 0, \eta_{ik}+\eta_{ij}\eta_{jk}\geq 0, \eta_{jk}+\eta_{ij}\eta_{ik}\geq 0
\end{aligned}
\end{equation}
for any face $\{ijk\}\in F$, then the generalized combinatorial curvature $K: B\rightarrow (0,+\infty)$ uniquely
determines the discrete conformal factor $u: B\rightarrow \mathbb{R}$.
\end{theorem}

\begin{remark}
The structure condition (\ref{structure condition}) is a direct consequence of the cosine law for generalized hyperbolic triangles with
all vertices hyper-ideal. Please refer to Section 5.1 in \cite{X2} and Remark \ref{expolation of structure condition}
in Section \ref{section 6} of this paper for the details of the geometric
explanation. The structure condition (\ref{structure condition}) has been previously
used in the study of discrete conformal structures on closed surfaces.
See, for instance, \cite{Z1,X1,X2a,X2} and others.
The rigidity result in Theorem \ref{main thm rigidity} can be taken to be
a complement of Guo-Luo's rigidity of the discrete conformal structures in \cite{GL}
and Guo's rigidity of vertex scaling in \cite{Guo2} on surface with boundary.
\end{remark}

\begin{figure}[!htb]
\centering
  \includegraphics[height=.55\textwidth,width=.5\textwidth]{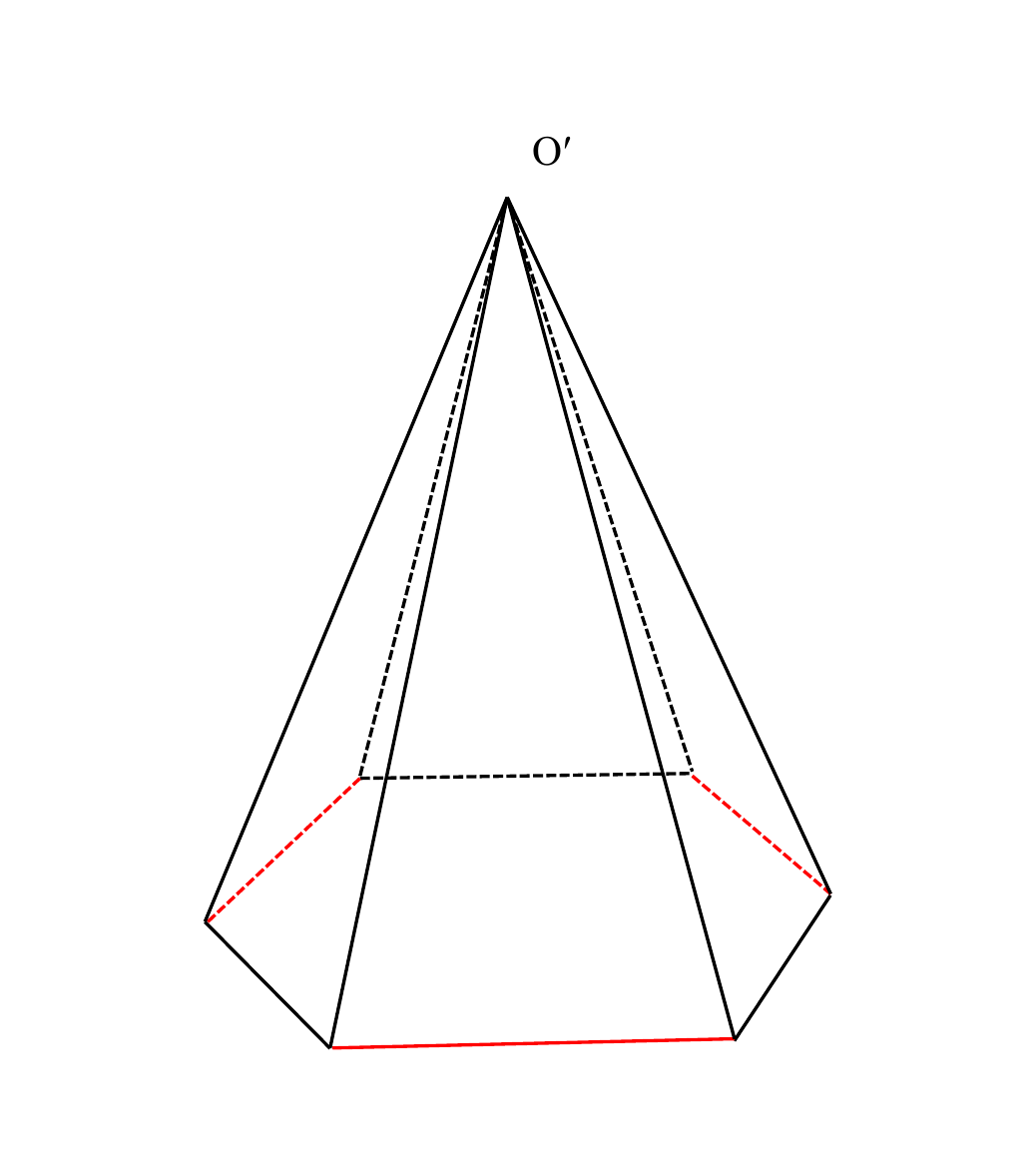}
  \caption{Generalized hyperbolic pyramid with a
right-angled hyperbolic hexagonal base}
  \label{pyramid}
\end{figure}

As a byproduct of Theorem \ref{main thm rigidity},
we prove some new convexity of the volume of generalized hyperbolic pyramids with
right-angled hyperbolic hexagonal bases.
Suppose  $Ov_iv_jv_k$ is a generalized hyperbolic tetrahedron satisfying (1) all the vertices $O$,$v_i$,$v_j$,$v_k$ are hyper-ideal,
(2) the intersection of the hyperbolic plane $P_O$ dual to $O$ and $\mathbb{H}^3$ is a generalized hyperbolic triangle
with all vertices hyper-ideal and edges intersecting with $\mathbb{H}^3$. Here and in the following, we use $P_v$ to denote the hyperbolic
plane dual to a hyper-ideal point $v$.
Please refer to Figure \ref{Generalized hyperbolic triangle for surfaces with boundary} for such a generalized hyperbolic tetrahedron.
Truncating the generalized  hyperbolic tetrahedron $Ov_iv_jv_k$ with the planes dual to $O$,$v_i$,$v_j$,$v_k$
gives rise to a generalized hyperbolic pyramid $C$, which has a base given by a right-angled hyperbolic hexagon and an apex $O'$.
Please refer to Figure \ref{pyramid} for the resulting generalized hyperbolic pyramid.
If $O'$ is hyper-ideal, we further truncate $C$ by the hyperbolic plane $P_{O'}$ dual to $O'$.
Otherwise, we keep the generalized hyperbolic pyramid $C$ invariant.
Denote the volume of the resulting generalized hyperbolic polyhedra by $V$ and
call it the volume of the generalized hyperbolic pyramid $C$.
Further denote the intersection angles of $P_O$
with $P_{v_i}$,$P_{v_j}$,$P_{v_k}$ by $\alpha_i$,  $\alpha_j$ and $\alpha_k$ respectively, which are dihedral angles
of the generalized hyperbolic pyramid $C$ at the edges $P_O\cap P_{v_i}$,$P_O\cap P_{v_j}$,$P_O\cap P_{v_k}$ respectively.

\begin{theorem}\label{main theorem convexity of volume}
The volume $V$ of the generalized hyperbolic pyramid $C$ is a strictly concave function of the
dihedral angles $\alpha_i$, $\alpha_j$ and $\alpha_k$.
\end{theorem}

Since Chow-Luo's pioneering work \cite{CL} on combinatorial Ricci flows for Thurston's circle packings
and Luo's work \cite{L1} on combinatorial Yamabe flow for vertex scaling of piecewise linear metrics on closed surfaces,
combinatorial curvature flows have been important approaches for
constructing geometrical structures on surfaces. There are lots of important works on combinatorial curvature
flows on surfaces. See, for instance, \cite{CL,L1,GLSW,GGLSW,X2, Ge-thesis,Ge,GH1,GX3,ZX,WX, Guo2,X3, LX} and others.
Aiming at finding hyperbolic metrics on surfaces with totally geodesic boundary components of prescribed lengths,
we introduce the following combinatorial curvature flows to deform the
discrete conformal structures in Definition \ref{defn of DCG on surf w bdy},
including combinatorial Ricci flow, combinatorial Calabi flow and fractional combinatorial Calabi flow.
Set
\begin{equation}\label{alpha parameter}
\alpha_i=\arctan e^{-u_i}, \alpha_i\in (0,\frac{\pi}{2}).
\end{equation}
By Definition \ref{defn of DCG on surf w bdy} and the formula (\ref{alpha parameter}),  we have
\begin{equation}\label{length of l_ij in alpha}
\cosh l_{ij}=\frac{\cos\alpha_i\cos \alpha_j+\eta_{ij}}{\sin\alpha_i\sin\alpha_j}.
\end{equation}
The formula (\ref{length of l_ij in alpha}) motivates a nice geometric interpolation
in $3$-dimensional hyperbolic geometry in Section \ref{section 6}
for the discrete conformal structures in Definition \ref{defn of DCG on surf w bdy}.
For simplicity, we also call the function $\alpha: B\rightarrow (0, \frac{\pi}{2})$ in (\ref{alpha parameter}) as a discrete conformal factor,
 if it causes no confusion in the context.
\begin{definition}
Suppose $(\Sigma, \mathcal{T})$ is an ideally triangulated surface with boundary
and $\eta: E\rightarrow (-1, +\infty)$ is a weight defined on the edges.
The combinatorial Ricci flow for the discrete conformal structures
in Definition \ref{defn of DCG on surf w bdy} on $(\Sigma, \mathcal{T}, \eta)$  is defined to be
    \begin{equation}\label{CRF}
\begin{aligned}
\left\{
  \begin{array}{ll}
    \frac{d\alpha_i}{dt}=\overline{K}_i -K_i,  &\hbox{ } \\
    \alpha(0)=\alpha_0,  &\hbox{ }
  \end{array}
\right.
\end{aligned}
\end{equation}
where $\overline{K}\in \mathbb{R}^n_{>0}$ is a positive function defined on $B=\{1, 2, \cdots, n\}$,
$\alpha_0$ is an admissible discrete conformal factor on $(\Sigma, \mathcal{T}, \eta)$.
The combinatorial Calabi flow for the discrete conformal structures in Definition \ref{defn of DCG on surf w bdy} on $(\Sigma, \mathcal{T}, \eta)$ is defined to be
\begin{equation}\label{CCF}
\begin{aligned}
\left\{
  \begin{array}{ll}
    \frac{d\alpha_i}{dt}=\Delta (K-\overline{K})_i, & \hbox{ } \\
    \alpha(0)=\alpha_0, & \hbox{ }
  \end{array}
\right.
\end{aligned}
\end{equation}
where $\Delta=-(\frac{\partial K}{\partial \alpha})$ is the discrete Laplace operator.
\end{definition}

In Proposition \ref{positive definiteness of global Jacobian},
we prove that the discrete Laplace operator $\Delta$ is strictly negative definite under the structure condition (\ref{structure condition}).
Following \cite{WX}, we define the fractional combinatorial Laplace operator $\Delta^s$ for $s\in \mathbb{R}$
as follows. Recall that a symmetric positive definite matrix $A_{n\times n}$ could be written as
$$A=P^T \cdot \text{diag}\{\lambda_1, \cdots, \lambda_n\}\cdot P,$$
where $P\in O(n)$ and $\lambda_1\leq  \cdots\leq \lambda_n$ are positive eigenvalues of $A$.
For any $s\in \mathbb{R}$, $A^s$ is defined to be
$$A^s=P^T \cdot \text{diag}\{\lambda_1^s, \cdots, \lambda_n^s\}\cdot P.$$
The $2s$-th order fractional discrete Laplace operator $\Delta^s$ is defined to be
\begin{equation}\label{fractional laplace}
  \Delta^s=-\left(\frac{\partial K}{\partial \alpha}\right)^s.
\end{equation}
Motivated by \cite{WX}, we introduce the following fractional combinatorial Calabi flow for the discrete conformal structures in Definition \ref{defn of DCG on surf w bdy} on $(\Sigma, \mathcal{T}, \eta)$
\begin{equation}\label{FCCF}
\begin{aligned}
\left\{
  \begin{array}{ll}
    \frac{d\alpha_i}{dt}=\Delta^s (K-\overline{K})_i, & \hbox{ } \\
    \alpha(0)=\alpha_0, & \hbox{ }
  \end{array}
\right.
\end{aligned}
\end{equation}
where $\Delta^s$ is the fractional discrete Laplace operator defined by (\ref{fractional laplace}).
If $s=0$, the fractional combinatorial Calabi flow (\ref{FCCF}) is reduced to the combinatorial Ricci flow (\ref{CRF}).
If $s=1$, the fractional combinatorial Calabi flow (\ref{FCCF}) is reduced to the combinatorial Calabi flow (\ref{CCF}).
We have the following result on the combinatorial curvatures flows.

\begin{theorem}\label{main theorem curvature flows}
Suppose $(\Sigma, \mathcal{T})$ is an ideally triangulated surface with boundary
and $\eta: E\rightarrow (-1, +\infty)$ is a weight defined on the edges satisfying the structure condition (\ref{structure condition}).
\begin{description}
  \item[(a)] The combinatorial Ricci flow (\ref{CRF}) and the combinatorial Calabi flow (\ref{CCF}) are negative gradient flows.
  \item[(b)] Suppose there exists an admissible discrete conformal structure $\overline{\alpha}$
  such that $K(\overline{\alpha})=\overline{K}$, then there exists a positive number $\delta$ such that
  if $||\alpha_0-\overline{\alpha}||<\delta$,  the solutions of
  the combinatorial Ricci flow (\ref{CRF}), the combinatorial Calabi flow (\ref{CCF}) and
  the fractional  combinatorial Calabi flow (\ref{FCCF})
  exist for all time and converge exponentially fast to $\overline{\alpha}$.
  \item[(c)] The solutions of  the combinatorial Ricci flow (\ref{CRF}), the combinatorial Calabi flow (\ref{CCF}) and
  the fractional  combinatorial Calabi flow (\ref{FCCF}) can not reach the boundary of the admissible space of discrete conformal
  factors $\alpha$ in $(0, \frac{\pi}{2})^n$.
\end{description}
\end{theorem}

The paper is organized as follows.
In Section \ref{section 2}, we give a characterization of the admissible space of discrete conformal factors on ideally
triangulated surfaces.
In Section \ref{section 3}, we prove that the Jacobian matrix of the generalized inner angles
with respect to the discrete conformal factors in a generalized hyperbolic triangle (right-angled hyperbolic hexagon)
is symmetric and positive definite.
 In Section \ref{section 4}, we prove the global rigidity of the  generalized combinatorial curvature, i.e. Theorem
 \ref{main thm rigidity}.
 In Section \ref{section 5}, we study the properties of the solutions to the combinatorial Ricci flow (\ref{CRF}),
 combinatorial Calabi flow (\ref{CCF}) and fractional combinatorial Calabi flow (\ref{FCCF})
on ideally triangulated surfaces with boundary and prove Theorem \ref{main theorem curvature flows}.
In Section \ref{section 6}, we study the relationships between
the discrete conformal structures in Definition \ref{defn of DCG on surf w bdy}
and $3$-dimensional hyperbolic space. As a result, we prove some
new convexity of the volume of some generalized hyperbolic pyramids with
right-angled hyperbolic hexagonal base in dihedral angles,
i.e. Theorem \ref{main theorem convexity of volume}.
In Section \ref{section 7}, we discuss some interesting open problems on hyperbolic discrete conformal structures on surfaces with boundary.
\\
\\
\textbf{Acknowledgements}\\[8pt]
The research of the author is supported by
the Fundamental Research Funds for the Central Universities under Grant no. 2042020kf0199.

\section{The admissible space of discrete conformal factors}\label{section 2}

We denote the space of functions $\alpha: B\rightarrow (0, \frac{\pi}{2})$ such that
\begin{equation}\label{inequ in alpha}
  \begin{aligned}
\frac{\cos\alpha_i\cos \alpha_j+\eta_{ij}}{\sin\alpha_i\sin\alpha_j}>1
  \end{aligned}
\end{equation}
for the edge $\{ij\}\in E$ as $\mathcal{W}^\alpha_{ij}$ and denote
the space of admissible discrete conformal factor $\alpha: B\rightarrow (0, \frac{\pi}{2})$ as $\mathcal{W}^\alpha$.

\begin{theorem}\label{thm admissible space}
Suppose $(\Sigma, \mathcal{T})$ is an ideally triangulated surface with boundary
and $\eta: E\rightarrow (-1, +\infty)$ is a weight defined on the edges.
For each edge $\{ij\}\in E$, the space $\mathcal{W}^\alpha_{ij}$
is a convex polytope in $(0,\frac{\pi}{2})^n$.
As a result, the admissible space
$$\mathcal{W}^\alpha=\cap_{\{ij\}\in E}\mathcal{W}^\alpha_{ij}$$
is a convex polytope in $(0,\frac{\pi}{2})^n$.
\end{theorem}
\proof
By the definition of $\alpha$ in (\ref{alpha parameter}), a function $\alpha: B\rightarrow (0, \frac{\pi}{2})$ belongs to $\mathcal{W}^\alpha_{ij}$
if and only if (\ref{inequ in alpha}) is valid,
which is equivalent to
\begin{equation}\label{equ proof admiss space}
  \cos(\alpha_i+\alpha_j)>-\eta_{ij}.
\end{equation}
If $-\eta_{ij}<-1$, i.e. $\eta_{ij}>1$,
the condition (\ref{equ proof admiss space}) is satisfied for any $\alpha_i, \alpha_j\in (0,\frac{\pi}{2})$.
If $-\eta_{ij}\geq-1$, i.e. $-1<\eta_{ij}\leq 1$, then the condition (\ref{equ proof admiss space})
implies that
$\alpha_i+\alpha_j<\arccos(-\eta_{ij}).$
In any case, for the edge $\{ij\}\in E$, the space
\begin{equation}\label{W_ij in alpha}
  \mathcal{W}^\alpha_{ij}=\{\alpha\in (0,\frac{\pi}{2})^n|\cos(\alpha_i+\alpha_j)>-\eta_{ij}\}
\end{equation}
is a convex polytope in $(0,\frac{\pi}{2})^2$.
\qed

\begin{remark}\label{remark on alpha domain}
The admissible space $\mathcal{W}^\alpha$ is nonempty. Especially, $\mathcal{W}^\alpha$ contains the points $\alpha\in (0, \frac{\pi}{2})^n$
with all $\alpha_i$ small enough.
One can also use (\ref{length of l_ij in alpha})
as the definition of discrete conformal structures on ideally triangulated surfaces with boundary with
$\alpha\in (0, \pi)$ as a discrete conformal factor instead of (\ref{DCS on bordered surface}), which corresponds to $\alpha\in (0, \frac{\pi}{2})$.
Following the proof of Theorem \ref{thm admissible space},
we can also prove that the admissible space of discrete conformal factors $\alpha$
is a convex polytope in $(0,\pi)^n$ in this case. This definition of discrete conformal structures is reasonable from the viewpoint
of 3-dimensional hyperbolic geometry in Section \ref{section 6}.
However, we can not prove the rigidity for the generalized combinatorial curvature in this setting.
Please refer to Remark \ref{remark on non-rigidity}.
\end{remark}


\section{Jacobian matrix of the generalized angles in a generalized hyperbolic triangle}\label{section 3}

Suppose $\{ijk\}\in F$ is a generalized hyperbolic triangle with hyper-ideal vertices $i,j,k$, which corresponds to a
right-angled hyperbolic hexagon adjacent the boundary components $i,j,k\in B$.
The length of the boundary arc in $\{ijk\}\in F$ facing the edge $\{jk\}$ is called the
generalized angle of the generalized hyperbolic triangle $\{ijk\}$ at $i$.
Denote the edge lengths of three nonadjacent edges $\{ij\}, \{ik\}, \{jk\}$ as $l_{ij}, l_{ik}, l_{jk}$ respectively
and the generalized inner angles at the hyper-ideal vertices $i,j,k$ as $\theta_i$, $\theta_j$, $\theta_k$ respectively.

\subsection{Symmetry of the Jacobian matrix}
Set
\begin{equation*}
  \gamma_i=\eta_{jk}+\eta_{ij}\eta_{ik}.
\end{equation*}
We have the following result on the symmetry of the Jacobian matrix $\frac{\partial (\theta_i, \theta_j, \theta_k)}{\partial (\alpha_i, \alpha_j, \alpha_k)}$

\begin{lemma}\label{symmetry of Jacobian}
For discrete conformal factor $\alpha\in \mathcal{W}^\alpha_{ij}\cap\mathcal{W}^\alpha_{ik}\cap\mathcal{W}^\alpha_{jk}$,
\begin{equation*}
  \begin{aligned}
\frac{\partial \theta_i}{\partial \alpha_j}
=\frac{\partial \theta_j}{\partial \alpha_i}
=\frac{1}{A_{ijk}\sinh^2 l_{ij}\sin^2\alpha_i\sin^2\alpha_j\sin\alpha_k}
 [(1-\eta^2_{ij})\cos\alpha_k+\gamma_i\cos\alpha_j+\gamma_j\cos\alpha_i],
  \end{aligned}
\end{equation*}
where $A_{ijk}=\sinh l_{ij}\sinh l_{ik}\sinh \theta_i$.
\end{lemma}
\proof
By the derivative cosine law for  right-angled hyperbolic  hexagons \cite{GL}, we have
\begin{equation}\label{equ 1 proof in symmetry}
  \begin{aligned}
\frac{\partial \theta_i}{\partial l_{ij}}=\frac{-\sinh l_{jk}\cosh \theta_j}{A_{ijk}},\
\frac{\partial \theta_i}{\partial l_{ik}}=\frac{-\sinh l_{jk}\cosh \theta_k}{A_{ijk}},\
\frac{\partial \theta_i}{\partial l_{jk}}=\frac{\sinh l_{jk}}{A_{ijk}}.
  \end{aligned}
\end{equation}
By the formula (\ref{length of l_ij in alpha}) of the hyperbolic length in discrete conformal factor $\alpha$, we have
\begin{equation}\label{equ 2 proof in symmetry}
  \begin{aligned}
\frac{\partial l_{ij}}{\partial \alpha_j}=-\frac{\cos\alpha_i+\eta_{ij}\cos\alpha_j}{\sinh l_{ij}\sin\alpha_i\sin^2\alpha_j},\
 \frac{\partial l_{ik}}{\partial \alpha_j}=0,\
 \frac{\partial l_{jk}}{\partial \alpha_j}=-\frac{\cos\alpha_k+\eta_{jk}\cos\alpha_j}{\sinh l_{jk}\sin^2\alpha_j\sin\alpha_k}.
  \end{aligned}
\end{equation}
By the chain rules, we have
\begin{equation}\label{equ 3 proof in symmetry}
  \begin{aligned}
\frac{\partial \theta_i}{\partial \alpha_j}
=\frac{\partial \theta_i}{\partial l_{ij}}\frac{\partial l_{ij}}{\partial \alpha_j}
+\frac{\partial \theta_i}{\partial l_{ik}} \frac{\partial l_{ik}}{\partial \alpha_j}
+\frac{\partial \theta_i}{\partial l_{jk}}\frac{\partial l_{jk}}{\partial \alpha_j}
=\frac{\partial \theta_i}{\partial l_{ij}}\frac{\partial l_{ij}}{\partial \alpha_j}
+\frac{\partial \theta_i}{\partial l_{jk}}\frac{\partial l_{jk}}{\partial \alpha_j}.
  \end{aligned}
\end{equation}
Submitting (\ref{equ 1 proof in symmetry}) and (\ref{equ 2 proof in symmetry}) into (\ref{equ 3 proof in symmetry}) gives
\begin{equation}\label{equ 4 proof in symmetry}
  \begin{aligned}
\frac{\partial \theta_i}{\partial \alpha_j}
=&\frac{\sinh l_{jk}\cosh \theta_j}{A_{ijk}}\cdot\frac{\cos\alpha_i+\eta_{ij}\cos\alpha_j}{\sinh l_{ij}\sin\alpha_i\sin^2\alpha_j}
-\frac{\sinh l_{jk}}{A_{ijk}}\cdot\frac{\cos\alpha_k+\eta_{jk}\cos\alpha_j}{\sinh l_{jk}\sin^2\alpha_j\sin\alpha_k}\\
=&\frac{\cos\alpha_i+\eta_{ij}\cos\alpha_j}{A_{ijk}\sinh^2 l_{ij}\sin\alpha_i\sin^2\alpha_j}[\cosh l_{ij}\cosh l_{jk}+\cosh l_{ik}]
  -\frac{\cos\alpha_k+\eta_{jk}\cos\alpha_j}{A_{ijk}\sin^2\alpha_j\sin\alpha_k},
  \end{aligned}
\end{equation}
where the cosine law for right-angled hyperbolic  hexagons is used in the last line.
Submitting (\ref{length of l_ij in alpha}) into (\ref{equ 4 proof in symmetry}), by lengthy but direct calculations,
we have
\begin{equation}\label{equ 5 proof in symmetry}
  \begin{aligned}
\frac{\partial \theta_i}{\partial \alpha_j}
=&\frac{1}{A_{ijk}\sinh^2 l_{ij}\sin^2\alpha_i\sin^2\alpha_j\sin\alpha_k}\\
 &\cdot [(1-\eta^2_{ij})\cos\alpha_k+(\eta_{ik}+\eta_{ij}\eta_{jk})\cos\alpha_i +(\eta_{jk}+\eta_{ij}\eta_{ik})\cos\alpha_j].
  \end{aligned}
\end{equation}
Note that $A_{ijk}$ is symmetric in $i,j,k$ by the sine law for right-angled hyperbolic  hexagons, we have
$\frac{\partial \theta_i}{\partial \alpha_j}=\frac{\partial \theta_j}{\partial \alpha_i}$ by  (\ref{equ 5 proof in symmetry}).
\qed

\begin{remark}
  Lemma \ref{symmetry of Jacobian} is still valid if we use   (\ref{length of l_ij in alpha})
as the definition of discrete conformal structures on ideally triangulated surfaces with boundary with $\alpha\in (0, \pi)$.
\end{remark}

As a direct corollary of Lemma \ref{symmetry of Jacobian}, we have the following result.
\begin{corollary}
If the weight
$\eta\in (-1,1]$ and satisfies the structure condition (\ref{structure condition}),
then
$$\frac{\partial \theta_i}{\partial \alpha_j}\geq 0,$$
the equality of which is attained if and only if $\eta_{ij}=1$ and $\eta_{ik}+\eta_{jk}=0$.
\end{corollary}


In the special case of $\eta\equiv 0$, we further have the following interesting formula on the relationships between
$\frac{\partial \theta_i}{\partial \alpha_i}$, $\frac{\partial \theta_i}{\partial \alpha_j}$ and
$\frac{\partial \theta_i}{\partial \alpha_k}$.

\begin{lemma}\label{lemma GT formula for eta 0}
If $\eta\equiv 0$, for discrete conformal factor $\alpha\in \mathcal{W}^\alpha_{ij}\cap\mathcal{W}^\alpha_{ik}\cap\mathcal{W}^\alpha_{jk}$, we have
\begin{equation}\label{GT 's formula for eta=0}
  \begin{aligned}
 \frac{\partial \theta_i}{\partial \alpha_i}
 = \frac{\partial \theta_i}{\partial \alpha_j}\cosh l_{ij}+\frac{\partial \theta_i}{\partial \alpha_k}\cosh l_{ik}.
  \end{aligned}
\end{equation}
\end{lemma}
\proof
By Lemma \ref{symmetry of Jacobian}, we have
\begin{equation}\label{equ 1 in relation of partial deriv}
  \begin{aligned}
\frac{\partial \theta_i}{\partial \alpha_j}
=\frac{\cot\alpha_k}{A_{ijk}\sinh^2 l_{ij}\sin^2\alpha_i\sin^2\alpha_j},
\frac{\partial \theta_i}{\partial \alpha_k}
=\frac{\cot\alpha_j}{A_{ijk}\sinh^2 l_{ik}\sin^2\alpha_i\sin^2\alpha_k}
  \end{aligned}
\end{equation}
in the case of $\eta\equiv 0$.
By the chain rules, we have
\begin{equation}\label{equ 2 in relation of partial deriv}
  \begin{aligned}
\frac{\partial \theta_i}{\partial \alpha_i}
=&\frac{\partial \theta_i}{\partial l_{ij}}\frac{\partial l_{ij}}{\partial \alpha_i}
+\frac{\partial \theta_i}{\partial l_{ik}} \frac{\partial l_{ik}}{\partial \alpha_i}
+\frac{\partial \theta_i}{\partial l_{jk}}\frac{\partial l_{jk}}{\partial \alpha_i}\\
=&\frac{\sinh l_{jk}\cosh \theta_j}{A_{ijk}}\frac{\cos\alpha_j+\eta_{ij}\cos\alpha_i}{\sinh l_{ij}\sin^2\alpha_i\sin\alpha_j}
    +\frac{\sinh l_{jk}\cosh \theta_k}{A_{ijk}}\frac{\cos\alpha_k+\eta_{ik}\cos\alpha_i}{\sinh l_{ik}\sin^2\alpha_i\sin\alpha_k}\\
=&\frac{1}{A_{ijk}\sinh^2 l_{ij}\sin^2\alpha_i\sin\alpha_j}(\cosh l_{ij}\cosh l_{jk}+\cosh l_{ik})(\cos\alpha_j+\eta_{ij}\cos\alpha_i)\\
   &+\frac{1}{A_{ijk}\sinh^2 l_{ik}\sin^2\alpha_i\sin\alpha_k}(\cosh l_{ik}\cosh l_{jk}+\cosh l_{ij})(\cos\alpha_k+\eta_{ik}\cos\alpha_i)
  \end{aligned}
\end{equation}
Submitting (\ref{length of l_ij in alpha}) and $\eta\equiv 0$ into (\ref{equ 2 in relation of partial deriv}), by lengthy but direct calculations,
we have
\begin{equation}\label{equ 3 in relation of partial deriv}
  \begin{aligned}
\frac{\partial \theta_i}{\partial \alpha_i}
=\frac{\cot \alpha_i\cot \alpha_j\cot \alpha_k}{A_{ijk}\sinh^2 l_{ij}\sin^2\alpha_i\sin^2\alpha_j}
   +\frac{\cot \alpha_i\cot \alpha_j\cot \alpha_k}{A_{ijk}\sinh^2 l_{ik}\sin^2\alpha_i\sin^2\alpha_k}.
  \end{aligned}
\end{equation}
Note that
\begin{equation}\label{equ 4 in relation of partial deriv}
  \begin{aligned}
\cosh l_{ij}=\cot \alpha_i\cot \alpha_j, \cosh l_{ik}=\cot \alpha_i\cot \alpha_k
  \end{aligned}
\end{equation}
in the case of $\eta\equiv 0$ by (\ref{length of l_ij in alpha}).
Combining the formulae  (\ref{equ 1 in relation of partial deriv}),  (\ref{equ 3 in relation of partial deriv}) and  (\ref{equ 4 in relation of partial deriv}) gives the formula (\ref{GT 's formula for eta=0}).
\qed

\begin{remark}
The formula (\ref{GT 's formula for eta=0}) was first obtained by Glickenstein-Thomas \cite{GT} (see also \cite{Thomas})
for generic hyperbolic discrete conformal structures on closed surfaces, which has lots of applications.
See, for instance, \cite{WX,X2, X2a, XZ1, XZ2} and others.
This is the first time the formula (\ref{GT 's formula for eta=0}) proved for hyperbolic discrete conformal structures on surfaces with boundary.
It is conceived that
the formula (\ref{GT 's formula for eta=0}) holds for any hyperbolic discrete conformal structures
on ideally triangulated surfaces with boundary.
\end{remark}

\subsection{Positive definiteness of the Jacobian matrix}
The aim of this subsection is to prove that the Jacobian matrix $\frac{\partial (\theta_i, \theta_j, \theta_k)}{\partial (u_i, u_j,u_k)}$
for a generalized triangle $\{ijk\}\in F$ is positive definite for admissible discrete conformal factors.
\begin{lemma}\label{nonsingular of Jacobian}
Suppose $\{ijk\}\in F$ is a generalized  triangle adjacent to $i,j,k\in B$ and
$\eta\in (-1, +\infty)$ is a weight defined on the edges satisfying the structure condition (\ref{structure condition}),
then the Jacobian matrix
$\frac{\partial (\theta_i, \theta_j, \theta_k)}{\partial (\alpha_i, \alpha_j,\alpha_k)}$
is nonsingular
with positive determinant
for any discrete conformal factor $\alpha\in \mathcal{W}^\alpha_{ij}\cap\mathcal{W}^\alpha_{ik}\cap\mathcal{W}^\alpha_{jk}$.
\end{lemma}
\proof
By the chain rules, we have
\begin{equation}\label{equ 0 proof nonsingular}
  \begin{aligned}
\frac{\partial (\theta_i, \theta_j, \theta_k)}{\partial (\alpha_i, \alpha_j,\alpha_k)}
=\frac{\partial (\theta_i, \theta_j, \theta_k)}{\partial (l_{ij}, l_{ik},l_{jk})}
\cdot\frac{\partial (l_{ij}, l_{ik},l_{jk})}{\partial (\alpha_i, \alpha_j,\alpha_k)}.
  \end{aligned}
\end{equation}
By the proof of Lemma \ref{symmetry of Jacobian},
we have
\begin{equation}\label{equ 1 proof nonsingular}
  \begin{aligned}
  \left(
    \begin{array}{c}
      d\theta_i \\
      d\theta_j  \\
      d\theta_k  \\
    \end{array}
  \right)
  =&\frac{-1}{\sinh l_{ij}\sinh l_{ik}\sinh \theta_i} \\
  &\cdot\left(
                                                                  \begin{array}{ccc}
                                                                    \sinh l_{jk} & 0 & 0 \\
                                                                    0 & \sinh l_{ik} & 0 \\
                                                                    0 & 0 & \sinh l_{ij} \\
                                                                  \end{array}
                                                                \right)
   \left(
      \begin{array}{ccc}
        \cosh \theta_j & \cosh \theta_k & -1 \\
        \cosh \theta_i  & -1 & \cosh \theta_k \\
       -1  & \cosh \theta_i & \cosh \theta_j \\
      \end{array}
    \right)
    \left(
      \begin{array}{c}
        dl_{ij} \\
        dl_{ik} \\
        dl_{jk} \\
      \end{array}
    \right),
  \end{aligned}
\end{equation}
which implies
\begin{equation*}
  \begin{aligned}
\det \left(\frac{\partial (\theta_i, \theta_j, \theta_k)}{\partial (l_{ij}, l_{ik},l_{jk})}\right)
=\left(\frac{-1}{\sinh l_{ij}\sinh l_{ik}\sinh \theta_i}\right)^3\sinh l_{ij}\sinh l_{ik}\sinh l_{jk} \det M
  \end{aligned}
\end{equation*}
with $M$ being the last $3\times 3$ matrix in (\ref{equ 1 proof nonsingular}).
By direct calculations, we have
\begin{equation*}
  \begin{aligned}
\det M
=-\cosh ^2\theta_i-\cosh ^2\theta_j-\cosh ^2\theta_k-2\cosh \theta_i\cosh \theta_j\cosh \theta_k+1<0,
  \end{aligned}
\end{equation*}
which implies that
\begin{equation}\label{equ 2 proof nonsingular}
  \begin{aligned}
  \det \left(\frac{\partial (\theta_i, \theta_j, \theta_k) }{\partial (l_{ij}, l_{ik}, l_{jk}) }\right)>0.
  \end{aligned}
\end{equation}
On the other hand, by the formula (\ref{length of l_ij in alpha}) of hyperbolic length in $\alpha$, we have
\begin{equation*}
  \begin{aligned}
\left(
  \begin{array}{c}
    dl_{ij} \\
    dl_{ik} \\
    dl_{jk} \\
  \end{array}
\right)
=-\left(
                                                                  \begin{array}{ccc}
                                                                    \frac{1}{\sinh l_{ij}} & 0 & 0 \\
                                                                    0 & \frac{1}{\sinh l_{ik}} & 0 \\
                                                                    0 & 0 & \frac{1}{\sinh l_{jk}} \\
                                                                  \end{array}
                                                                \right)
\cdot J \cdot
 \left(
   \begin{array}{c}
     d\alpha_i \\
     d\alpha_j \\
     d\alpha_k \\
   \end{array}
 \right),
  \end{aligned}
\end{equation*}
where $J$ is the following $3\times 3$ matrix
\begin{equation*}
  \begin{aligned}
J=\left(
   \begin{array}{ccc}
     \cosh l_{ij}\cot \alpha_i+\cot \alpha_j & \cosh l_{ij}\cot \alpha_j+\cot \alpha_i & 0 \\
     \cosh l_{ik}\cot \alpha_i+\cot \alpha_k & 0 & \cosh l_{ik}\cot \alpha_k+\cot \alpha_i \\
     0 & \cosh l_{jk}\cot \alpha_j+\cot \alpha_k & \cosh l_{jk}\cot \alpha_k+\cot \alpha_j \\
   \end{array}
 \right).
  \end{aligned}
\end{equation*}
By lengthy but direct calculations, we have
\begin{equation}\label{equ 3 proof nonsingular}
  \begin{aligned}
&\det \left(\frac{\partial (l_{ij}, l_{ik}, l_{jk})}{\partial (\alpha_i, \alpha_j,\alpha_k)}\right)\\
=&\frac{-1}{\sinh l_{ij}\sinh l_{ik}\sinh l_{jk}}\det J\\
=&\frac{1}{\sinh l_{ij}\sinh l_{ik}\sinh l_{jk}\sin^3\alpha_i\sin^3\alpha_j\sin^3\alpha_j}\\
 &\cdot [2(1+\eta_{ij}\eta_{ik}\eta_{jk})\cos \alpha_i\cos\alpha_j\cos \alpha_k
 +\gamma_i\cos \alpha_i(\cos^2 \alpha_j+\cos^2 \alpha_k)\\
  &\ \ \ +\gamma_j\cos \alpha_j(\cos^2 \alpha_i+\cos^2 \alpha_k) +\gamma_k\cos \alpha_k(\cos^2 \alpha_i+\cos^2 \alpha_j)].
  \end{aligned}
\end{equation}
Note that $\alpha_i, \alpha_j, \alpha_k\in (0,\frac{\pi}{2})$ by (\ref{alpha parameter})
and $\gamma_i\geq 0, \gamma_j\geq 0,  \gamma_k\geq 0$ by the structure condition (\ref{structure condition}),
we have
\begin{equation}\label{equ 4 proof nonsingular}
  \begin{aligned}
&\det \left(\frac{\partial (l_{ij}, l_{ik}, l_{jk})}{\partial (\alpha_i, \alpha_j,\alpha_k)}\right)\\
\geq&\frac{2\cos\alpha_i\cos\alpha_j\cos\alpha_k}{\sinh l_{ij}\sinh l_{ik}\sinh l_{jk}\sin^3\alpha_i\sin^3\alpha_j\sin^3\alpha_j}\\
&\cdot (1+\eta_{ij}\eta_{ik}\eta_{jk}+\eta_{jk}+\eta_{ij}\eta_{ik}+\eta_{ij}+\eta_{ik}\eta_{jk}+\eta_{ik}+\eta_{ij}\eta_{jk})\\
=&\frac{2\cos\alpha_i\cos\alpha_j\cos\alpha_k}{\sinh l_{ij}\sinh l_{ik}\sinh l_{jk}\sin^3\alpha_i\sin^3\alpha_j\sin^3\alpha_j}(1+\eta_{ij})(1+\eta_{ik})(1+\eta_{jk})\\
>&0
  \end{aligned}
\end{equation}
by (\ref{equ 3 proof nonsingular}),
where the conditions $\eta\in (-1, +\infty)$ and $\alpha_i, \alpha_j, \alpha_k\in (0,\frac{\pi}{2})$ are used in the last line.
Combining the equations (\ref{equ 2 proof nonsingular}) and (\ref{equ 4 proof nonsingular}) gives
$$\det \left(\frac{\partial (\theta_i, \theta_j, \theta_k)}{\partial (\alpha_i, \alpha_j,\alpha_k)}\right)>0$$
by (\ref{equ 0 proof nonsingular}), which implies that the Jacobian matrix
$\frac{\partial (\theta_i, \theta_j, \theta_k)}{\partial (\alpha_i, \alpha_j,\alpha_k)}$
is nonsingular.
\qed

\begin{remark}\label{remark on non-rigidity}
By Remark \ref{remark on alpha domain},
we can define the discrete conformal structure using the formula (\ref{length of l_ij in alpha})
with $\alpha_i\in (0, \pi)$. However, in this case, we do not have $\cos \alpha_i>0$, which is technically
necessary in the proof of Lemma \ref{nonsingular of Jacobian}.
\end{remark}


To prove that the Jacobian matrix $\frac{\partial (\theta_i, \theta_j, \theta_k)}{\partial (\alpha_i, \alpha_j,\alpha_k)}$
is positive definite, following \cite{X2,X2a}, we further introduce the following parameterized admissible space
\begin{equation*}
  \mathcal{A}_{ijk}=\{(\alpha_i, \alpha_j, \alpha_k, \eta_{ij}, \eta_{ik}, \eta_{jk})\in \mathcal{W}^\alpha_{ijk}(\eta)\times (-1, +\infty)^3|
\gamma_i\geq 0, \gamma_j\geq 0, \gamma_k\geq 0\}
\end{equation*}
for a generalized triangle $\{ijk\}\in F$,
where we use $\mathcal{W}^\alpha_{ijk}(\eta)$ to denote
$\mathcal{W}^\alpha_{ij}\cap\mathcal{W}^\alpha_{ik}\cap\mathcal{W}^\alpha_{jk}$ depending on $\eta$
for simplicity.
The parameterized admissible space $\mathcal{A}_{ijk}$ can be taken as a fibre bundle over the space
$$\Gamma:=\{(\eta_{ij}, \eta_{ik}, \eta_{jk})\in (-1, +\infty)^3|\gamma_i\geq 0, \gamma_j\geq 0, \gamma_k\geq 0\}$$
with the fibre given by $\mathcal{W}^\alpha_{ijk}(\eta)$ over $\eta\in \Gamma$.
\begin{lemma}[\cite{X2a} Lemma 2.7]\label{connectness of Gamma}
  The space $\Gamma$ is path connected.
\end{lemma}

As a direct corollary of Lemma \ref{connectness of Gamma}, we have the following result.

\begin{corollary}\label{connectness of para adm space}
Suppose $\{ijk\}\in F$ is a generalized triangle adjacent to $i,j,k\in B$.
Then the parameterized admissible space $\mathcal{A}_{ijk}$ is connected.
\end{corollary}
\proof
Set
$$f_{ij}(\alpha, \eta)=\frac{\cos\alpha_i\cos \alpha_j+\eta_{ij}}{\sin\alpha_i\sin\alpha_j}.$$
Then
$(\alpha_i, \alpha_j, \alpha_k, \eta_{ij}, \eta_{ik}, \eta_{jk})\in \mathcal{A}_{ijk}$
if and only if $f_{ij}>1$,  $f_{ik}>1$, $f_{jk}>1$.
By the continuity of $f_{ij}, f_{ik}, f_{jk}$,
if $(\alpha_i, \alpha_j, \alpha_k, \eta_{ij}, \eta_{ik}, \eta_{jk})\in \mathcal{A}_{ijk}$,
then there is a convex neighborhood $U$ of $(\alpha_i, \alpha_j, \alpha_k, \eta_{ij}, \eta_{ik}, \eta_{jk})$
such that $U\subseteq\mathcal{A}_{ijk}$.
As a result, for any fixed point $(\overline{\alpha}_i, \overline{\alpha}_j, \overline{\alpha}_k)\in \mathcal{W}^\alpha_{ijk}(\eta_0)$,
there is a connected neighborhood $V$ of $\eta_0$ in $\Gamma$
such that the space
\begin{equation*}
  \{(\overline{\alpha}_i, \overline{\alpha}_j, \overline{\alpha}_k, \eta_{ij}, \eta_{ik}, \eta_{jk})\in \mathcal{A}_{ijk} |
  (\overline{\alpha}_i, \overline{\alpha}_j, \overline{\alpha}_k)\in \mathcal{W}^\alpha_{ijk}(\eta), (\eta_{ij}, \eta_{ik}, \eta_{jk})\in V\}
\end{equation*}
is connected.
Then the connectivity of the parameterized admissible space $\mathcal{A}_{ijk}$
follows from Theorem \ref{thm admissible space} and Lemma \ref{connectness of Gamma}.
\qed


\begin{theorem}\label{positive definiteness of Jacobian}
Suppose $\{ijk\}\in F$ is a generalized hyperbolic triangle adjacent to $i,j,k\in B$
and $\eta\in (-1, +\infty)$ is a weight defined on the edges satisfying the structure condition (\ref{structure condition}).
Then the Jacobian matrix
  $$\frac{\partial (\theta_i, \theta_j,\theta_k)}{\partial (\alpha_i, \alpha_j, \alpha_k)}$$
is strictly positive definite for any discrete conformal factor in
$\mathcal{W}^\alpha_{ij}\cap\mathcal{W}^\alpha_{ik}\cap\mathcal{W}^\alpha_{jk}$.
\end{theorem}
\proof
Note that the Jacobian matrix $\frac{\partial (\theta_i, \theta_j,\theta_k)}{\partial (\alpha_i, \alpha_j, \alpha_k)}$
can be taken as a matrix-valued function defined on the parameterized admissible space $\mathcal{A}_{ijk}$.
By Lemma \ref{nonsingular of Jacobian} and Corollary \ref{connectness of para adm space},
the eigenvalues of $\frac{\partial (\theta_i, \theta_j,\theta_k)}{\partial (\alpha_i, \alpha_j, \alpha_k)}$
are nonzero continuous functions defined on the connected parameterized admissible space $\mathcal{A}_{ijk}$.
To prove that $\frac{\partial (\theta_i, \theta_j,\theta_k)}{\partial (\alpha_i, \alpha_j, \alpha_k)}$
is positive definite, we just need to choose a point $p$ in $\mathcal{A}_{ijk}$ such that
$\frac{\partial (\theta_i, \theta_j,\theta_k)}{\partial (\alpha_i, \alpha_j, \alpha_k)}$
is positive definite at $p$.

Take $p=(\alpha_i, \alpha_j, \alpha_k, 0, 0, 0)\in \mathcal{A}_{ijk}$.
By Lemma \ref{symmetry of Jacobian}, we have
\begin{equation}\label{equ in proof of posive definite}
  \frac{\partial \theta_i}{\partial \alpha_j}> 0, \frac{\partial \theta_i}{\partial \alpha_k}>0
\end{equation}
at $p$.
By Lemma \ref{lemma GT formula for eta 0}, we further have
$$ \frac{\partial \theta_i}{\partial \alpha_i}
 > \frac{\partial \theta_i}{\partial \alpha_j}+\frac{\partial \theta_i}{\partial \alpha_k}>0$$
at $p$ by (\ref{equ in proof of posive definite}), which implies that
the Jacobian matrix
$\frac{\partial (\theta_i, \theta_j,\theta_k)}{\partial (\alpha_i, \alpha_j, \alpha_k)}$
is diagonal dominant and then positive definite at $p$.
\qed

\section{Rigidity of discrete conformal structures}\label{section 4}
In this section, we give a proof of Theorem \ref{main thm rigidity}.
\bigskip
\\
\textbf{Proof of Theorem \ref{main thm rigidity}:}
By Lemma \ref{symmetry of Jacobian},
$\theta_id\alpha_i+\theta_jd\alpha_j+\theta_kd\alpha_k$
is a smooth closed  $1$-form on $\mathcal{W}^\alpha_{ij}\cap\mathcal{W}^\alpha_{ik}\cap\mathcal{W}^\alpha_{jk}$.
By Theorem \ref{thm admissible space}, the function
\begin{equation*}
  \mathcal{E}_{ijk}(\alpha_i, \alpha_j,\alpha_k)=\int^{(\alpha_i, \alpha_j,\alpha_k)}\theta_id\alpha_i+\theta_jd\alpha_j+\theta_kd\alpha_k
\end{equation*}
is a well-defined smooth function on $\mathcal{W}^\alpha_{ij}\cap\mathcal{W}^\alpha_{ik}\cap\mathcal{W}^\alpha_{jk}$
with $\nabla \mathcal{E}_{ijk}=(\theta_i, \theta_j, \theta_k)$,
which is strictly convex by Theorem \ref{positive definiteness of Jacobian}.
Set
\begin{equation}\label{energy function on surface}
 \mathcal{E}(\alpha)=\sum_{\{ijk\}\in F}\mathcal{E}_{ijk}(\alpha_i, \alpha_j,\alpha_k)
\end{equation}
for $\alpha\in \mathcal{W}^\alpha$.
Then $\mathcal{E}(\alpha)$ is a strictly convex smooth function defined on the convex admissible space $\mathcal{W}^\alpha$
with the gradient given by
$$\nabla\mathcal{E}=K.$$
Then the global rigidity of the generalized combinatorial curvature $K$ follows from the following well-known result in analysis.
\begin{lemma}\label{injectivity of the gradient of convex function}
If $W: \Omega\rightarrow \mathbb{R}$ is a $C^2$-smooth strictly convex function defined on a convex domain $\Omega\subseteq \mathbb{R}^n$,
then its gradient $\nabla W: \Omega\rightarrow \mathbb{R}^n$ is injective.
\end{lemma}
\qed

In the proof of Theorem \ref{main thm rigidity},
we have proved the following result on the Jacobian matrix $(\frac{\partial K}{\partial \alpha})$.

\begin{proposition}\label{positive definiteness of global Jacobian}
Suppose $(\Sigma, \mathcal{T})$ is an ideally triangulated surface with boundary
and $\eta: E\rightarrow (-1, +\infty)$ is a weight defined on the edges satisfying the structure condition (\ref{structure condition}).
Then the Jacobian matrix $(\frac{\partial K}{\partial \alpha})$
is symmetric and strictly positive definite on the admissible space $\mathcal{W}^\alpha$.
\end{proposition}

\section{Combinatorial curvature flows on surfaces with boundary}\label{section 5}
In this section, we study some basic properties of
the combinatorial Ricci flow (\ref{CRF}), the combinatorial Calabi flow (\ref{CCF}), the fractional combinatorial Calabi flow (\ref{FCCF})
and give a proof of Theorem \ref{main theorem curvature flows}.

\begin{lemma}\label{CRF is a gradient flow}
  The combinatorial Ricci flow (\ref{CRF}) is a negative gradient flow of the convex energy function
  defined by
  \begin{equation*}
  \overline{\mathcal{E}}(\alpha)=\mathcal{E}(\alpha)-\sum_{i\in B}\overline{K}_i\alpha_i,
\end{equation*}
where $\mathcal{E}(\alpha)$ is defined by the formula (\ref{energy function on surface}).
\end{lemma}
\proof
By the proof of Theorem \ref{main thm rigidity},
$\overline{\mathcal{E}}(\alpha)$ is a smooth convex function with
$$\nabla \overline{\mathcal{E}}=K-\overline{K},$$
which implies that the combinatorial Ricci flow (\ref{CRF}) is
a negative gradient flow of the convex energy function $\overline{\mathcal{E}}(\alpha)$.
\qed

\begin{corollary}\label{mononicity of energy along CRF}
The energy function $\overline{\mathcal{E}}(\alpha)$ is decreasing along the combinatorial Ricci flow (\ref{CRF}).
Furthermore, the generalized combinatorial Calabi energy $\mathcal{C}(\alpha)$ defined by
\begin{equation*}
  \mathcal{C}(\alpha)=\frac{1}{2}\sum_{i\in B}(K_i-\overline{K}_i)^2
\end{equation*}
is decreasing along the  combinatorial Ricci flow (\ref{CRF}).
\end{corollary}
\proof
The monotonicity of $\overline{\mathcal{E}}(\alpha)$ along the combinatorial Ricci flow (\ref{CRF}) follows from
Lemma \ref{CRF is a gradient flow}.
For the generalized combinatorial Calabi energy $\mathcal{C}(\alpha)$,
by direct calculations, we have
$$\frac{d\mathcal{C}(\alpha(t))}{dt}=\sum_{i\in B}(K_i-\overline{K}_i)\frac{dK_i}{dt}
=-(K-\overline{K})^T\cdot \left(\frac{\partial K}{\partial \alpha}\right)\cdot (K-\overline{K})\leq 0,$$
along the  combinatorial Ricci flow (\ref{CRF}) by Proposition \ref{positive definiteness of global Jacobian},
which is strictly negative unless $K=\overline{K}$.
\qed

\begin{lemma}\label{CCF is a gradient flow}
  The combinatorial Calabi flow (\ref{CRF})
  is a negative gradient flow of the generalized combinatorial Calabi energy energy $\mathcal{C}(\alpha)$.
As a result, the generalized combinatorial Calabi energy $\mathcal{C}(\alpha)$
is decreasing along the  combinatorial Calabi flow (\ref{CCF}).
Furthermore, the energy function $\overline{\mathcal{E}}(\alpha)$
is decreasing along the  combinatorial Ricci flow (\ref{CRF}).
\end{lemma}
\proof
By direct calculations, we have
$$\nabla_{\alpha_i}\mathcal{C}=\sum_{j\in B}\frac{\partial K_j}{\partial \alpha_i}(K_j-\overline{K}_j)=-\Delta (K-\overline{K})_i$$
by Proposition \ref{positive definiteness of global Jacobian},
which implies that the combinatorial Calabi flow (\ref{CRF}) is a negative gradient flow of $\mathcal{C}(w)$.
Similarly, we have
$$\frac{d\overline{\mathcal{E}}(\alpha(t))}{dt}=\sum_{i\in B}\nabla_{\alpha_i} \overline{\mathcal{E}}\cdot \frac{d\alpha_i}{dt}
=-(K-\overline{K})^T\cdot \left(\frac{\partial K}{\partial \alpha}\right)\cdot (K-\overline{K})\leq 0$$
by Proposition \ref{positive definiteness of global Jacobian},
which is strictly negative unless $K=\overline{K}$.
\qed

Lemma \ref{CRF is a gradient flow} and Lemma \ref{CCF is a gradient flow}
proves Theorem \ref{main theorem curvature flows} (a).
Following the arguments in the proof of Lemma \ref{CRF is a gradient flow}, Corollary \ref{mononicity of energy along CRF}
and Lemma \ref{CCF is a gradient flow}, we have the following result on the fractional combinatorial Calabi flow (\ref{FCCF}).

\begin{lemma}\label{mononicity of FCCF}
Suppose $(\Sigma, \mathcal{T})$ is an ideally triangulated surface with boundary
and $\eta: E\rightarrow (-1, +\infty)$ is a weight defined on the edges satisfying the structure condition (\ref{structure condition}).
Then for any $s\in \mathbb{R}$, the energy function $\overline{\mathcal{E}}(\alpha)$ and
the generalized combinatorial Calabi energy $\mathcal{C}(\alpha)$
is decreasing along the  fractional combinatorial Calabi flow (\ref{FCCF}).
\end{lemma}

\begin{remark}
For generic $s\in \mathbb{R}$, except $s=0, 1$, the fractional combinatorial Calabi flow (\ref{FCCF}) is not a gradient flow.
Furthermore, as the fractional discrete Laplace operator $\Delta^s$ is generically a non-local operator,
the definition of which involves the eigenvalues of matrices,
the fractional combinatorial Calabi flow (\ref{FCCF}) is generically (except $s=0, 1$) a non-local combinatorial curvature flow.
\end{remark}

Now we give a proof of Theorem \ref{main theorem curvature flows} (b).
As the proofs for the combinatorial Ricci flow (\ref{CRF}), the combinatorial Calabi flow (\ref{CCF}) and
  the fractional combinatorial Calabi flow (\ref{FCCF}) are similar,
  we only give the proof for the fractional combinatorial Calabi flow (\ref{FCCF}) for simplicity.

\textbf{Proof of Theorem \ref{main theorem curvature flows} (b):}
Set $\Gamma(\alpha)=\Delta^s (K-\overline{K})$ for the fractional combinatorial Calabi flow (\ref{FCCF}).
Then $\overline{\alpha}$ is an equilibrium point of the system (\ref{FCCF}) by assumption
and
$$D\Gamma|_{\alpha=\overline{\alpha}}=-\left(\frac{\partial K}{\partial \alpha}\right)(\overline{\alpha}),$$
which is strictly negative definite by Proposition \ref{positive definiteness of global Jacobian}.
This implies that $\overline{\alpha}$ is a local attractor of the fractional combinatorial Calabi flow (\ref{FCCF}).
Then the long time existence of the solution $\alpha(t)$ to (\ref{FCCF}) and
the exponential convergence of the solution $\alpha(t)$ to $\overline{\alpha}$
follows from the Lyapunov stability theorem (\cite{P}, Chapter 5).
\qed


Recall the characterization (\ref{W_ij in alpha}) of the space $\mathcal{W}^\alpha_{ij}$
of discrete conformal factors $\alpha: B\rightarrow (0, \frac{\pi}{2})$
such that (\ref{inequ in alpha}) is satisfied for an edge $\{ij\}\in E$.
In the case of $\eta\in (-1, 1]$, $\mathcal{W}^\alpha_{ij}$ has a non-empty boundary $\partial\mathcal{W}^\alpha_{ij}$
in $(0, \frac{\pi}{2})^n$ defined by
\begin{equation*}
  \partial\mathcal{W}^\alpha_{ij}=\{\alpha\in (0,\frac{\pi}{2})^n|\alpha_i+\alpha_j=\arccos(-\eta_{ij})\}.
\end{equation*}

\begin{lemma}\label{estimate for generalized combinatorial curvature at an edge}
Assume $(\Sigma, \mathcal{T})$ is an ideally triangulated surface with boundary and
$\eta: E\rightarrow (-1, +\infty)$ is a weight on the edges with $\eta_{ij}\in (-1, 1]$ for some edge $\{ij\}\in E$.
  For any $M>0$, there exists a positive constant $\epsilon_{ij}=\epsilon_{ij}(M)$ such that
if $\alpha\in \mathcal{W}^\alpha$ satisfies
$$\alpha_i+\alpha_j< \arccos(-\eta_{ij})+\epsilon_{ij},$$
then the generalized combinatorial curvature $K$ satisfies
$$K_i>M, K_j>M.$$
\end{lemma}
\proof
Suppose $\{ijk\}\in F$ is a face adjacent to the edge $\{ij\}\in E$.
By the cosine law for right-angled hyperbolic  hexagons, we have
\begin{equation*}
\begin{aligned}
\cosh \theta_i
=\frac{\cosh l_{ij}\cosh l_{ik}+\cosh l_{jk}}{\sinh l_{ij}\sinh l_{ik}}
>\frac{\cosh l_{ij}\cosh l_{ik}}{\sinh l_{ij}\sinh l_{ik}}
>\frac{\cosh l_{ij}}{\sinh l_{ij}},
\end{aligned}
\end{equation*}
which implies that
$\theta_i\rightarrow+\infty$ uniformly as $l_{ij}\rightarrow 0^+$.
Note that $K_i\geq\theta_i$ by the definition of generalized combinatorial curvature of discrete hyperbolic metrics
on ideally triangulated surfaces with boundary
 and
$l_{ij}\rightarrow 0^+$ is equivalent to $\alpha\in \mathcal{W}^\alpha$ and
$\alpha_i+\alpha_j\rightarrow (\arccos(-\eta_{ij}))^-$,
we have $K_i\rightarrow+\infty$ uniformly as
$\alpha_i+\alpha_j\rightarrow (\arccos(-\eta_{ij}))^-$.
The same arguments show that $K_j\rightarrow+\infty$ uniformly as $\alpha_i+\alpha_j\rightarrow (\arccos(-\eta_{ij}))^-$.
Therefore, for any number $M>0$, there exists a positive constant $\epsilon_{ij}=\epsilon_{ij}(M)$ such that
if $\alpha\in \mathcal{W}^\alpha$ satisfies $\alpha_i+\alpha_j< \arccos(-\eta_{ij})+\epsilon_{ij}$,
then $K_i>M, K_j>M.$
\qed

As an application of Lemma \ref{estimate for generalized combinatorial curvature at an edge}, we have the following result,
which is equivalent to Theorem \ref{main theorem curvature flows} (c).
\begin{proposition}\label{alpha(t) can not in W espsilon}
Assume $(\Sigma, \mathcal{T})$ is an ideally triangulated surface with boundary and
$\eta: E\rightarrow (-1, +\infty)$ is a weight on the edges with $\eta_{ij}\in (-1, 1]$ for some edges $\{ij\}\in E$.
Let  $\bar{K}\in \mathbb{R}^n_{>0}$ be a function defined on the boundary components $B$.
For any number $s\in \mathbb{R}$ and any initial value $\alpha_0\in \mathcal{W}^\alpha$,
there exists a constant $\epsilon = \epsilon(s, \alpha_0, \bar{K}) > 0$
such that the solution $\alpha(t)$ to the fractional combinatorial Calabi flow (\ref{FCCF}) can never be in the region
$$\mathcal{W}^\alpha_\epsilon = \{\alpha\in \mathcal{W}^\alpha|d(\alpha, \partial \mathcal{W})< \epsilon\},$$
where $d$ is the standard Euclidean metric on $\mathbb{R}^n$.
\end{proposition}
\proof
The proof is paralleling to that of Lemma 2.8 in \cite{LX}. For completeness, we give the proof here.
Set
$$M = \max_{i\in B} \{|\bar{K}_i|+ \sqrt{2\mathcal{C}(\alpha_0)}\}.$$
Suppose $\eta_{ij}\in (-1,1]$ for some edge $\{ij\}\in E$.
By Lemma \ref{estimate for generalized combinatorial curvature at an edge},
there exists $\epsilon_{ij}=\epsilon_{ij}(M)>0$ such that if
$$\alpha_i+\alpha_j < \arccos(-\eta_{ij})+2\epsilon_{ij},$$
then
 $$K_i(\alpha)>M, K_j(\alpha)>M.$$
Set
$$\epsilon_0=\min_{\{ij\}\in E, \eta_{ij}\in (-1,1]} \epsilon_{ij}>0.$$
Then if $\alpha\in \mathcal{W}^\alpha$ satisfies
$$\alpha_i+\alpha_j < \arccos(-\eta_{ij}) +2 \epsilon_0$$
for some edge $\{ij\}\in E$ with $\eta_{ij}\in (-1,1]$, we have
$K_i(\alpha)>M,$
which further implies that
\begin{equation}\label{equ in proof of solu in W alpha}
  |K_i(\alpha) - \bar{K}_i| \geq |K_i(\alpha)| - |\bar{K}_i| > M - |\bar{K}_i|\geq\sqrt{2\mathcal{C}(\alpha_0)}.
\end{equation}

We claim that the solution $\alpha(t)$ to the fractional combinatorial Calabi flow (\ref{FCCF})  can never be in the region $\mathcal{W}^\alpha_{\epsilon_0}$.
Otherwise, there exists some $t_0\in [0, +\infty)$ and an edge $\{ij\} \in E$ with $\eta_{ij}\in (-1,1]$ such that the solution $\alpha(t)$
to the fractional combinatorial Calabi flow (\ref{FCCF})
satisfies $\alpha(t_0)\in \mathcal{W}^\alpha$ and
$$\alpha_i(t_0)+\alpha_j(t_0) <  \arccos(-\eta_{ij}) +2 \epsilon_0,$$
which further implies that
\begin{equation}\label{equ 1 in proof of compactnees}
  |K_i(\alpha(t_0))- \bar{K}_i|>\sqrt{2\mathcal{C}(\alpha_0)}
\end{equation}
by (\ref{equ in proof of solu in W alpha}).
Note that the generalized combinatorial Calabi energy $\mathcal{C}(\alpha)$
is decreasing along the fractional combinatorial Calabi flow (\ref{FCCF})
by Lemma \ref{mononicity of FCCF}.
Therefore,  for any $t>0$, the solution $\alpha(t)$ to the fractional combinatorial Calabi flow (\ref{FCCF})  satisfies
$$|K_i(t) - \bar{K}_i| \leq \sqrt{2\mathcal{C}(\alpha(t))} \leq \sqrt{2\mathcal{C}(\alpha_0)}$$
for any $i\in B$, which contradicts (\ref{equ 1 in proof of compactnees}).
Therefore, the solution $\alpha(t)$ to the fractional combinatorial Calabi flow (\ref{FCCF})
can never be in the region $\mathcal{W}^\alpha_{\epsilon_0}$.
\qed

\begin{remark}
  The result in Proposition \ref{alpha(t) can not in W espsilon}
  is independent of the assumption on the existence of $\bar{\alpha}\in \mathcal{W}^\alpha$
with $K(\bar{\alpha})=\bar{K}$ in Theorem \ref{main theorem curvature flows}.
\end{remark}

\section{Relationships with 3-dimensional hyperbolic geometry}\label{section 6}

\subsection{Construction of generalized hyperbolic triangles}
\begin{figure}[!htb]
\begin{minipage}[t]{0.5\linewidth}
\centering
  \includegraphics[height=6.3cm,width=5cm]{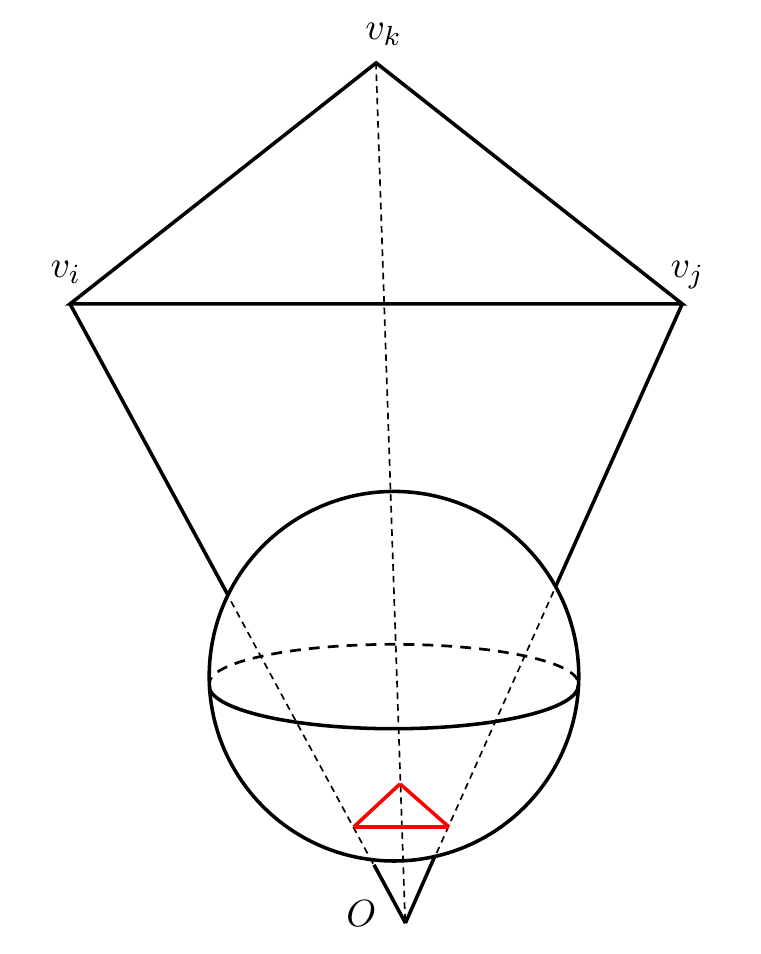}
  \caption{Hyperbolic triangle}
  \label{Hyperbolic triangle for closed surfaces}
\end{minipage}
\hfill
\begin{minipage}[t]{0.5\linewidth}
\centering
  \includegraphics[height=6.6cm,width=5.8cm]{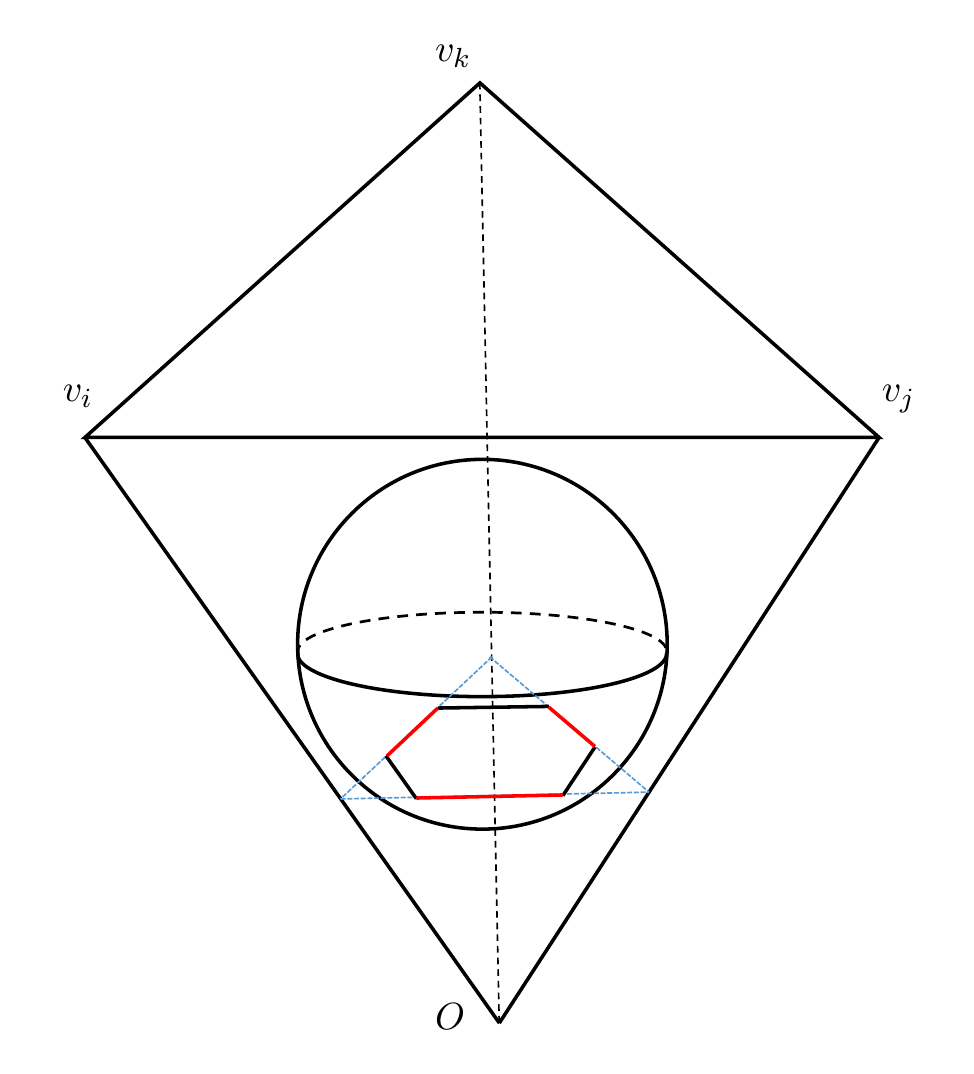}
  \caption{Generalized hyperbolic triangle}
  \label{Generalized hyperbolic triangle for surfaces with boundary}
\end{minipage}
\end{figure}

The key to define a discrete conformal structures on surfaces is to construct a (generalized) geometric triangle
with variables defined on the vertices and prescribed weights defined on the edges of a topological triangle, which is
closely related to $3$-dimensional hyperbolic geometry.
The relationships between discrete conformal structures on closed surfaces and $3$-dimensional hyperbolic geometry
were first observed by Bobenko-Pinkall-Springborn \cite{BPS} in the case of Luo's vertex scaling of piecewise linear metrics.
Let us give a quick review of Bobenko-Pinkall-Springborn's observation.
Suppose $Ov_iv_jv_k$ is an ideal hyperbolic tetrahedron in $\mathbb{H}^3$ with each ideal vertex attached with a horosphere,
which is usually referred as a decorated ideal hyperbolic tetrahedron.
Bobenko-Pinkall-Springborn found that Luo's constuction of Euclidean triangle via vertex scaling
corresponds exactly to the Euclidean triangle given by the intersection of $Ov_iv_jv_k$
and the horosphere at the ideal vertex $O$,
if the generalized edge lengths of the decorated ideal hyperbolic tetrahedron $Ov_iv_jv_k$ are properly assigned.
Based on this observation, Bobenko-Pinkall-Springborn \cite{BPS}
further introduced the vertex scaling for piecewise hyperbolic metrics
by perturbing the ideal vertex $O$ of the ideal hyperbolic tetrahedron $Ov_iv_jv_k$ to hyper-ideal while keeping the other vertices ideal.
In this case, the hyperbolic triangle is given by the intersection
of the generalized hyperbolic tetrahedron $Ov_iv_jv_k$ and the hyperbolic plane $P_O$ dual to the hyper-ideal vertex $O$
with the generalized lengths of the generalized hyperbolic tetrahedron $Ov_iv_jv_k$ properly assigned.
Please refer to Figure \ref{Hyperbolic triangle for closed surfaces} with $v_i, v_j, v_k$ ideal for the construction of hyperbolic vertex scaling.
The readers are suggested to refer to
Bobenko-Pinkall-Springborn's work \cite{BPS} for more details on this.

Motivated by Bobenko-Pinkall-Springborn's observations \cite{BPS},
Zhang-Guo-Zeng-Luo-Yau-Gu \cite{ZGZLYG} further constructed all $18$ types of
discrete conformal structures on closed surfaces in different background geometries by
perturbing the ideal vertices of the ideal hyperbolic tetrahedron $Ov_iv_jv_k$ to be hyperbolic, ideal or hyper-ideal.
Specially, for the hyperbolic discrete conformal structures on closed surfaces,
the vertex $O$ is required to be hyper-ideal and the lines $Ov_i$, $Ov_j$, $Ov_k$ are required to \textbf{intersect} with
the $3$-dimensional hyperbolic space $\mathbb{H}^3$. The constructed hyperbolic triangle is then given by the intersection of
the generalized hyperbolic tetrahedron  $Ov_iv_jv_k$ and the hyperbolic plane $P_O$ dual to $O$.
Please refer to Figure \ref{Hyperbolic triangle for closed surfaces} for the hyperbolic triangle in the Klein model constructed in this approach.

The discrete conformal structures on surfaces with boundary in Definition \ref{defn of DCG on surf w bdy}
are constructed by further perturbing the vertices of the generalized tetrahedron $Ov_iv_jv_k$ as follows.
Note that in Zhang-Guo-Zeng-Luo-Yau-Gu's construction of hyperbolic triangles, the vertex $O$ is hyper-ideal and
the lines $Ov_i$, $Ov_j$, $Ov_k$ are required to intersect with the $3$-dimensional hyperbolic space $\mathbb{H}^3$.
If we further perturb the vertex $O$ such that the lines $Ov_i$, $Ov_j$, $Ov_k$
do \textbf{NOT intersect} with $\mathbb{H}^3\cup \partial \mathbb{H}^3$ and
the intersection of the hyperbolic plane $P_O$ dual to $O$ with $Ov_iv_jv_k$ is a generalized hyperbolic triangle
with all vertices hyper-ideal and edges intersecting with $\mathbb{H}^3$,
then the intersection of $P_O$ with $Ov_iv_jv_k$ is exactly the generalized hyperbolic triangle
induced by a right-angled hyperbolic hexagon shown in Figure \ref{Colored right-angled hyperbolic  hexagon}.
Note that in this case, all the vertices of the generalized hyperbolic tetrahedron $Ov_iv_jv_k$ are hyper-ideal.
Please refer to Figure \ref{Generalized hyperbolic triangle for surfaces with boundary} for the construction.

  \begin{figure}[!htb]
\centering
  \includegraphics[height=.32\textwidth,width=1\textwidth]{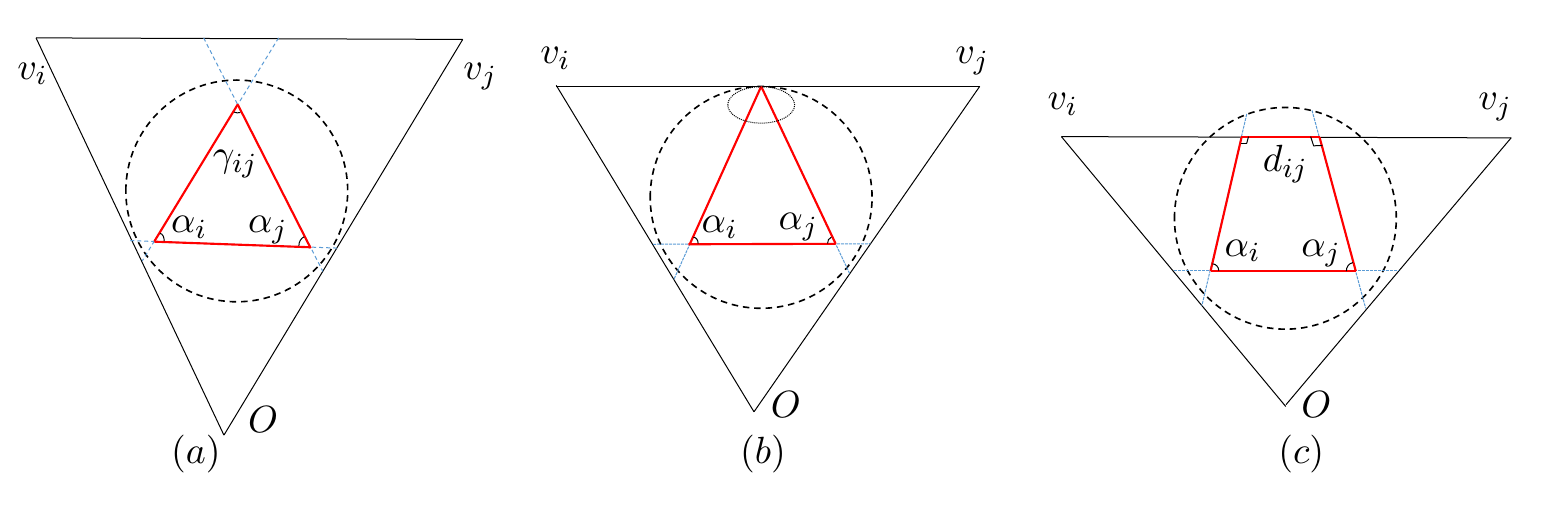}
  \caption{Lateral generalized triangles}
  \label{Lateral generalized triangles}
\end{figure}

Now we derive the formula (\ref{length of l_ij in alpha}) in the definition of the discrete conformal structure in
Definition \ref{defn of DCG on surf w bdy} using the construction above.
For this, we just need to consider a lateral generalized triangle $Ov_iv_j$ of the generalized hyperbolic tetrahedron $Ov_iv_jv_k$.
Denote the hyperbolic lines dual to $O, v_i, v_j$ as $L_O, L_i, L_j$ respectively.
By the requirement that $Ov_i$, $Ov_j$ do not intersect with $\mathbb{H}^2\cup \partial \mathbb{H}^2$, we can suppose
that $L_O$ intersects $L_i, L_j$ with the angles $\alpha_i, \alpha_j$ respectively. Please refer to Figure \ref{Lateral generalized triangles}.
If the line $v_iv_j$ does not intersect with $\mathbb{H}^2\cup \partial \mathbb{H}^2$, we can suppose $L_i$ and $L_j$ intersect with angle $\gamma_{ij}$. Please refer to Figure \ref{Lateral generalized triangles} (a) for this. In this case, we have
\begin{equation}\label{equ 1 in definition}
  \cosh l_{ij}=\frac{\cos\alpha_i\cos \alpha_j+\cos\gamma_{ij}}{\sin\alpha_i\sin\alpha_j}
\end{equation}
by the hyperbolic cosine law for hyperbolic triangles.
If the line $v_iv_j$ is tangential to $\partial\mathbb{H}^2$, then $L_O, L_i, L_j$ forms a generalized hyperbolic triangle with one ideal vertex and two hyperbolic vertices. Please refer to Figure \ref{Lateral generalized triangles} (b) for this.
In this case, we have
\begin{equation}\label{equ 2 in definition}
\cosh l_{ij}=\frac{\cos\alpha_i\cos \alpha_j+1}{\sin\alpha_i\sin\alpha_j}
\end{equation}
by the cosine law for generalized hyperbolic triangles with one ideal vertex and two hyperbolic vertices.
If the line $v_iv_j$ intersects with $\mathbb{H}^2$, then $L_O, L_i, L_j$ forms a generalized hyperbolic triangle with one hyper-ideal vertex and two hyperbolic vertices. Furthermore, there is a unique hyperbolic segment $L_{ij}$ perpendicular to $L_i$ and $L_j$,
the length of which is denoted by $d_{ij}$.
Please refer to Figure \ref{Lateral generalized triangles} (c) for this.
In this case, we have
\begin{equation}\label{equ 3 in definition}
\cosh l_{ij}=\frac{\cos\alpha_i\cos \alpha_j+\cosh d_{ij}}{\sin\alpha_i\sin\alpha_j}.
\end{equation}
The formulas (\ref{equ 1 in definition}), (\ref{equ 2 in definition}) and (\ref{equ 3 in definition})
together motivate us to define the hyperbolic length using the formula (\ref{length of l_ij in alpha}).
\begin{remark}\label{expolation of structure condition}
In the above construction, the weight $\eta_{ij}$ in the formula (\ref{length of l_ij in alpha})
is determined by the relative position of the two hyper-ideal vertices $v_i, v_j$.
Note that $v_i, v_j, v_k$ forms a generalized hyperbolic triangle with the vertices all hyper-ideal.
Following the arguments in Section 5.1 of \cite{X2}, the structure condition is a natural
consequence of the hyperbolic cosine law for such generalized hyperbolic triangles.
\end{remark}

\begin{remark}
 In \cite{GL}, Guo-Luo introduced some other types of discrete conformal structures on surfaces with boundary
  using Andreev-Thurston's circle packing approach with the standard hyperbolic cosine law
  replaced by different types of cosine laws in hyperbolic geometry.
  Guo-Luo's construction can be obtained by further perturbing
   the generalized hyperbolic tetrahedron $Ov_iv_jv_k$ in Figure \ref{Generalized hyperbolic triangle for surfaces with boundary}
  as follows.
  First, we truncated the generalized hyperbolic tetrahedron $Ov_iv_jv_k$ by the half space determined by $P_O$ containing $O$,
  which gives rise to a generalized hyperbolic polytope $v_iv_jv_kv_i'v_j'v_k'$ with $v_s'$ being the intersection of $Ov_s$ with $P_O$ for $s\in \{i,j,k\}$.
  Second, we further perturb the edges $v_iv_i'$, $v_jv_j'$, $v_kv_k'$ of the generalized hyperbolic polytope $v_iv_jv_kv_i'v_j'v_k'$
   such that $v_i, v_j,v_k$ becomes one point $O'$, $v_i', v_j',v_k'$ are kept
  hyper-ideal and $O'v_i'$, $O'v_j'$, $O'v_k'$ intersect with the hyperbolic space $\mathbb{H}^3$.
  Then Guo-Luo's definition of discrete conformal structures on surfaces with boundary
  corresponds to the generalized hyperbolic triangle $v_i'v_j'v_k'$ with
  the generalized length of $v_i'v_j'$ defined by the generalized length of $O'v_i'$, $O'v_j'$
  and the generalized angle $\angle v_i'O'v_j'$ using hyperbolic cosine laws.
One can refer to \cite{GL} for more details on Guo-Luo's construction.
From the arguments above, we can see that the discrete conformal structures in Definition \ref{defn of DCG on surf w bdy}
are dual to Guo-Luo's discrete conformal structures in \cite{GL}.
However, we do not know how to include Guo's vertex scaling of discrete hyperbolic metrics in \cite{Guo2}
using such geometric constructions.
\end{remark}

\subsection{Convexity of the volume of some generalized hyperbolic pyramids}

Suppose that $Ov_iv_jv_k$  is a generalized hyperbolic tetrahedron constructed above
for the discrete conformal structures in Definition \ref{defn of DCG on surf w bdy}
with the weights $\eta_{ij}, \eta_{ik},\eta_{jk}$ fixed.
$Ov_iv_jv_k$ can be attached with a generalized hyperbolic polytope in $\mathbb{H}^3$ as follows.
In the first step, we truncate the generalized
hyperbolic tetrahedron $Ov_iv_jv_k$ by $P_O, P_i, P_j, P_k$,
which gives rise to a generalized hyperbolic pyramid $C$ with a right-angled hyperbolic hexagonal base and an apex $O'$.
Please refer to Figure \ref{pyramid} for the generalized hyperbolic pyramid $C$.
If the resulting generalized hyperbolic pyramid $C$ still contains any hyper-ideal vertex,
then we continue the procedure in the first step until there is no hyper-ideal
vertex for the resulting hyperbolic polytope $P_f$.
Note that the final hyperbolic polytope $P_f$ may contain an ideal vertice at $O'$.
We define the volume of the generalized
hyperbolic pyramid $C$ to be the volume $V$ of  hyperbolic polytope $P_f$,
which is a function of $\alpha_i, \alpha_j, \alpha_k$.
By the generalized Schl\"{a}fli formula in \cite{R1} and the condition that the weights $\eta_{ij}, \eta_{ik},\eta_{jk}$ are fixed, we have
\begin{equation*}
  dV=-\frac{1}{2}(\theta_id\alpha_i+\theta_jd\alpha_j+\theta_kd\alpha_k),
\end{equation*}
which shows that the volume $V$ of the generalized hyperbolic pyramid $C$
is a strictly concave function of the parameters $\alpha_i, \alpha_j, \alpha_k$ by Theorem \ref{positive definiteness of Jacobian}.
In summary, we have the following result, which is equivalent to Theorem \ref{main theorem convexity of volume}.
\begin{proposition}\label{convexity of volume}
  Suppose $Ov_iv_jv_k$ is a generalized hyperbolic tetrahedron
  constructed for the discrete conformal structures in Definition \ref{defn of DCG on surf w bdy}
  with the relative positions of the hyper-ideal vertices $v_i, v_j, v_k$ fixed, i.e. the weights $\eta_{ij}, \eta_{ik},\eta_{jk}$ are  fixed.
  Then the volume $V$ of the generalized
hyperbolic pyramid $C$ constructed above
   is a strictly concave function of the dihedral angles $\alpha_i, \alpha_j, \alpha_k$.
\end{proposition}

\section{Open problems}\label{section 7}
\subsection{Classification of discrete conformal structures on surfaces with boundary}
There are different types of discrete conformal structures on closed surfaces that have been extensively studied in the history,
including tangential circle packings, Thurston's circle packings, inversive distance circle packings, Luo's vertex scaling,
generic discrete conformal structures proposed by Glickenstein et al. and others.
These different types of discrete conformal structures on closed surfaces are introduced and studied individually for a long time
until the works of Glickenstein \cite{G3}, Glickenstein-Thomas \cite{GT} and Zhang-Guo-Zeng-Luo-Yao-Gu \cite{ZGZLYG},
which unifies and generalizes different types of discrete conformal structures on closed surfaces.
Zhang-Guo-Zeng-Luo-Yao-Gu's approach \cite{ZGZLYG}
is motivated by Bobenko-Pinkall-Springborn's observations \cite{BPS} on the relationships
between Luo's vertex scaling of piecewise linear metrics and the $3$-dimensional hyperbolic geometry.
They explicitly constructed 18 different types of discrete conformal structures on closed surfaces by perturbing
the ideal vertices of ideal hyperbolic tetrahedron in the extended $3$-dimensional hyperbolic space.
Glickenstein and Glickenstein-Thomas's approach in \cite{G3, GT} is much different.
They defined the discrete conformal structures on closed surfaces by some reasonable axioms
and clarified the  discrete conformal structures they defined.
It is fantastical that the two different approaches give rise to the same discrete conformal structures on closed surfaces.
The rigidity of generic discrete conformal structures on closed surfaces was recently proved by the author in \cite{X2},
where the deformation of the discrete conformal structures was also studied.

Following Andreev-Thurston's approach,
Guo-Luo \cite{GL} introduced some other types of discrete conformal structures on surfaces with boundary
with the standard hyperbolic cosine law replaced by different types of cosine laws in hyperbolic geometry.
Following Luo's vertex scaling of piecewise linear metrics on closed surfaces,
Guo \cite{Guo2} also introduced the following discrete conformal structures on ideally triangulated surfaces with boundary,
called vertex scaling as well.

\begin{definition}[Guo \cite{Guo2}]\label{defn Guo's vertex scaling}
Suppose $(\Sigma, \mathcal{T})$ is an ideally triangulated surface with boundary.
Let $l$ and $l^0$ be two discrete hyperbolic metrics on $(\Sigma, \mathcal{T})$.
If there exists a function $u: B\rightarrow \mathbb{R}$ such that
\begin{equation*}
  \cosh \frac{l_{ij}}{2}=e^{u_i+u_j}\cosh \frac{l^0_{ij}}{2},
\end{equation*}
then the discrete hyperbolic metric $l$ is called vertex scaling of $l^0$. The function $u: B\rightarrow \mathbb{R}$ is called  a discrete conformal factor.
\end{definition}

Guo \cite{Guo2} proved the global rigidity of the vertex scaling in Definition \ref{defn Guo's vertex scaling} and studied the longtime behavior of the corresponding combinatorial Yamabe flow. See also \cite{LX, X3}.
Note that Guo's vertex scaling of discrete hyperbolic metrics on ideally triangulated surfaces with boundary in Definition \ref{defn Guo's vertex scaling} is formally different from the discrete conformal structure introduced in Definition \ref{defn of DCG on surf w bdy}.
As the discrete conformal structures on closed surfaces have been classified and their rigidities have been unified,
natural questions for discrete conformal structures on surfaces with boundary are as follows.

\begin{question}
  Can we find the full list of hyperbolic discrete conformal structures on ideally triangulated surfaces with boundary?
  Can we classify the hyperbolic discrete conformal structures on ideally triangulated surfaces with boundary
  following Glickenstein and Glickenstein-Thomas's axiomatic approach in \cite{G3, GT}?
Do the hyperbolic discrete conformal structures on ideally triangulated surfaces with boundary
have a unified version of rigidity as that in \cite{X2}?
\end{question}

\subsection{Prescribing the generalized combinatorial curvature on surfaces with boundary}

For discrete conformal structures on closed surfaces,
the prescribing combinatorial curvature problem has nice solutions.
For Thurston's circle packing, the image of the combinatorial curvature is a convex polytope, which
has been proved in Thurston's famous lecture notes \cite{T}.
For the vertex scaling, the prescribing combinatorial curvature problem have been perfectly solved by
Gu-Luo-Sun-Wu \cite{GLSW} in the Euclidean background geometry and by Gu-Guo-Luo-Sun-Wu \cite{GGLSW}
in the hyperbolic background geometry via introducing a new definition of discrete conformality allowing
the triangulations to be changed under the Delaunay condition. See also \cite{Sp} for the case of sphere.

Guo's vertex scaling on surfaces with boundary in Definition \ref{defn Guo's vertex scaling} is an analogue of Luo's vertex scaling
on closed surfaces.
Comparing Lemma \ref{symmetry of Jacobian} with Lemma 3.6 in \cite{X1}, one can see that the discrete conformal
structure on surfaces with boundary in Definition \ref{defn of DCG on surf w bdy}
is an analogue of the circle packings on closed surfaces.

A natural question related to prescribing generalized combinatorial curvature problem
is as follows.

\begin{question}
  Can we introduce a new definition of discrete conformality for discrete hyperbolic metrics on surfaces with boundary, following
  Gu-Luo-Sun-Wu \cite{GLSW} and Gu-Guo-Luo-Sun-Wu \cite{GGLSW}, and give a solution of the prescribing
  generalized combinatorial curvature problem?
\end{question}

Note that the prescribing generalized combinatorial curvature problem on surfaces with boundary
is equivalent to find a hyperbolic metric on surfaces with totally geodesic boundary components of prescribed lengths.
It is conceived that Luo's work \cite{L2} on Teichm\"{u}ller spaces of surfaces with boundary
will play a key role in the process.

\subsection{Long time behavior of the combinatorial curvature flows on surfaces with boundary}

The combinatorial Ricci (Yamabe) flow and combinatorial Calabi flow have been extensively studied on closed surfaces.
The corresponding long time existence and global convergence of the combinatorial curvature flows
has been well-established.
See, for instance, \cite{CL,L1,GLSW,GGLSW,X2} and others for the combinatorial Ricci (Yamabe) flow
and \cite{Ge-thesis,Ge,GH1,GX3,ZX,WX} and others for the combinatorial Calabi flow on closed surfaces.

For Guo's vertex scaling of discrete hyperbolic metrics on surfaces with boundary, the long time existence
and global convergence of combinatorial Yamabe flow is established in \cite{Guo2,X3}
and the long time existence
and global convergence of combinatorial Calabi flow is established in \cite{LX}.
However, for the hyperbolic discrete conformal structure on surfaces
with boundary in Definition \ref{defn of DCG on surf w bdy},
we only have the local convergence in Theorem \ref{main theorem curvature flows}.
A natural question related to the combinatorial curvature flows for the discrete conformal structures in Definition \ref{defn of DCG on surf w bdy}
is as follows.

\begin{question}
  Can we introduce some notion of surgery by flipping following \cite{GLSW, GGLSW} and
  prove the long time existence and global convergence of the combinatorial Ricci flow (\ref{CRF}),
  the combinatorial Calabi flow (\ref{CCF}) and the fractional combinatorial Calabi flow (\ref{FCCF})
  with surgery on surfaces with boundary?
\end{question}

(Xu Xu) School of Mathematics and Statistics, Wuhan University, Wuhan 430072, P.R. China

E-mail: xuxu2@whu.edu.cn\\[2pt]

\end{document}